\documentstyle[12pt]{article}

\def\qed{\hfill $\Box$}

\def\cz{ I \!\!\!\! C}

\def\gz{ Z \!\!\! Z}
\def\pr{ I \!\! P}
\def\t{\times}
\def\r{\rightarrow}

\def\s{\subset}
\def\se{\setminus}

\def\D{U}

\def\be{\beta}

\def\h\def\ph{\varphi}

\def\h{{\cal O}}

\def\noin{\noindent}
\def\be{\begin{equation}}
\def\ee{\end{equation}}
\def\vs{\bigskip}
\def\ms{\medskip}
\def\ss{\smallskip}
\begin{document}

\title {Plane curves with hyperbolic and C--hyperbolic complements}

\vspace{.5cm}

\author{ G. Dethloff, M. Zaidenberg}

\date{}
\maketitle

\vs

\noin We find sufficient conditions for the complement $\pr^2 \setminus C$ of a
plane curve $C$ to be C--hyperbolic. The latter means that some covering over
$\pr^2 \setminus C$ is Carath\'eodory hyperbolic. This  implies that this
complement  $\pr^2 \setminus C$ is Kobayashi hyperbolic, and (due to Lin's
Theorem) the fundamental group  $\pi_1 (\pr^2 \setminus C)$ does not contain a
nilpotent subgroup of finite index. We also give explicit examples of
irreducible such curves of any even degree $d \ge 6$.

\vs

\section{Introduction}

A complex space $X$ is said to be {\it C-hyperbolic} if there exists a
non-ramified covering $Y \r X$ such that $Y$ is Carath\'eodory hyperbolic,
i.e. the points in $Y$ are separated by bounded holomorphic functions (see
Kobayashi [21]).
If there exists a covering $Y$ of $X$ such that
for any point $p \in Y$ there exist only finitely many points $q \in Y$
which cannot be separated from $p$ by bounded holomorphic fuctions on $Y$,
then we say
that $X$ is  {\it almost C-hyperbolic}. There is a general problem:
{\it Which quasiprojective varieties are uniformized by bounded domains in
$\cz^n$?} In particular, such a variety must be C-hyperbolic. Here we
study plane projective curves whose complements are C-hyperbolic.
We prove the following

\vs

\noin {\bf 1.1. Theorem.} {\it Let $C \s \pr^2$ be an irreducible curve
of geometric genus $g$. Assume
that its dual curve $C^*$ is an immersed curve of degree $n$.

\noin a) If $g \ge 1$, then $\pr^2 \se C$ is C--hyperbolic.

\noin b) If $g = 0, \,n \ge 5$ and $C^*$ is a generic rational nodal curve,
then $\pr^2 \se C$ is almost C--hyperbolic.

\noin c) In both cases $\pr^2 \se C$ is Kobayashi complete hyperbolic
and hyperbolically embedded into $\pr^2$. }

\vs

Consider, for instance, an elliptic sextic with $9$ cusps (see (6.7)).
Such a sextic  can be given explicitly by
Schl\"afli's equation (see Gelfand, Kapranov and Zelevinsky [14]). It is dual to a smooth cubic and hence,
due to (a) its  complement is
C--hyperbolic. Actually, $6$ is the least possible degree of an
irreducible plane curve with C--hyperbolic complement (see (7.5)).

Note that C--hyperbolicity implies Kobayashi
hyperbolicity. S. Kobayashi [21] proposed the following

\vs

\noin {\bf Conjecture.} {\it Let ${\cal H}(d)$ be the set of all
hypersurfaces $D$ of degree $d$ in $\pr^n$ such that $\pr^n \se D$ is
complete hyperbolic and hyperbolically embedded into  $\pr^n$. Then for
any $d \geq 2n+1$ the set ${\cal H}(d)$ contains a Zariski open subset of
$\pr^{N(d)}$, where $\pr^{N(d)}$ is the complete linear system of effective
divisors of degree $d$ in $\pr^n$.} \vs

For $n > 2$ the problem is still open. For $n = 2$  Y.-T. Siu and
S.-K. Yeung [32] have announced a proof of the above conjecture in
particular case
of plane curves
of sufficiently large degree ($d > 10^6$). However, even for $n=2$ and for
small $d$ it is not so easy to construct explicit examples of
irreducible plane curves in ${\cal H}(d)$  (see Zaidenberg [37] and literature therein).
For the case of reducible curves see e.g. Dethloff, Schumacher and Wong [8,9].

The first examples of smooth curves in ${\cal H}(d)$ of any even degree
$d \geq 30$ were constructed by K. Azukawa and M. Suzuki [3]. A. Nadel [24]
mentioned such examples for any $d \geq 18$ which is divisible by 6. K. Masuda
and J. Noguchi [23] obtained smooth curves in ${\cal H}(d)$ for any $d \geq
21$.

In Zaidenberg [36] the existence of smooth curves in ${\cal H}(d)$ is proven (by
deformation arguments) for arbitrary $d \geq 5$; however, their equations are
not quite explicit. For instance, the equation of a smooth quintic in ${\cal
H}(5)$
includes five parameters which should be chosen successively small enough,
with unexplicit upper bounds.

In a series of papers by M. Green [16], J. Carlson and M. Green [4] and
H. Grauert and U. Peternell [15] sufficient conditions were found for
irreducible plane curves of genus $g \ge 2$ to be in ${\cal H}(d)$.
This leads to examples of irreducible (but singular) curves in ${\cal H}(d)$
with $d \ge 9$ (see the remark after (6.5)).

Generalizing the method of Green [16] (see the proof of Theorem 3.1, a)), we obtain
for any even $d \geq 6$
families of irreducible curves in ${\cal H}(d)$ described in terms
of genus and  singularities. While in all the examples known before the
curves were
of genus at least two, now we obtain such examples of elliptic or
rational curves.
They are all singular, and the method used is not
available to get such examples of smooth curves, even of higher genus.
On the other hand, it is clear that
an elliptic or a rational curve with hyperbolic complement must be singular.

We have presented above a family of elliptic sextics with C-hyperbolic
and hyperbolically embedded complements. Another example of curves with
hyperbolically embedded complements in
degree $6$, is the family of rational sextics with four nodes and
six cusps, where the cusps are on a conic (6.11). Such a curve is dual to
a generic rational nodal quartic,
and therefore, one can easily write down its explicit equation.

\vs

In fact,  C--hyperbolicity is a much stronger property than
Kobayashi hyperbolicity. To show this, recall that {\it a Liouville}
complex space is a space $Y$ such that all the bounded holomorphic
functions on $Y$ are constant. This property is just opposite of being
Carath\'eodory hyperbolic. By Lin's Theorem (see Lin [22], Theorem B), a Galois
covering $Y$ of a quasiprojective variety $X$ is Liouville if its group of
deck transformations is {\it almost nilpotent}, i.e. it contains a nilpotent
subgroup of finite index. It follows that, as soon as the fundamental
group $\pi_1( \pr^2 \se C)$ is
almost nilpotent, any covering over $\pr^2 \se C$ is
a Liouville one. In particular, this is so for a nodal (not necessarily
irreducible) plane curve  $C$ \footnote{i.e. a curve  $C$ with only normal
crossing singularities.}. Indeed, due to the Deligne-Fulton Theorem
(see Deligne [7] and Fulton [13]), in the latter case the group $\pi_1(\pr^2 \se C)$ is abelian.
As a corollary we obtain that for curves mentioned in Theorem 1.1 the
group $\pi_1( \pr^2 \se C)$ is not almost nilpotent (see Proposition 7.1).
For $n = {\rm deg}\, C^* \ge 2g +1$ this group even contains a free
subgroup with two generators (see sect. 7.a).

This shows that C--hyperbolicity of $\pr^2 \se C$ can be easily
destroyed under small deformations of $C$, by passing to a smooth or nodal
approximating curve $C'$. Observe that hyperbolicity of projective complements
is often stable and, in particular, smooth curves with hyperbolic complements
form an open subset (see Zaidenberg [36]). Whereas the locus of curves with C--hyperbolic
complements is contained in the locus of curves with singularities worse
than ordinary double points.

\ms

The paper is organized as follows. In section 2 we summarize necessary
background on plane algebraic curves, hyperbolic complex analysis and
on PGL$(2, \cz)$--actions on $\pr^n$. In section 3 we formulate Theorem 3.1
which is a generalization of Theorem 1.1. Its proof is given in sections
3--5. In this theorem we give sufficient conditions of C-hyperbolicity of
a complement of a plane curve together with its {\it artifacts}, i.e.
certain of its inflectional and cuspidal tangent lines. A curve has no
artifacts exactly when its dual is an immersed curve.

Section 6 is devoted to examples of curves of low
degrees with hyperbolic and C--hyperbolic complements. In section 7.a we
discuss the fundamental groups of the complements of curves with immersed dual.
In Proposition 7.5 we prove that $6$ is the minimal degree of irreducible
curves with C--hyperbolic complement. We also establish
genericity of inflexional tangents (i.e. artifacts) of a generic plane curve
(Proposition 7.6).

\ss

A part of the results of this paper was reported at the Hayama Conference on
Geometric Complex Analysis (Japan, March 1995; see Dethloff and Zaidenberg  [10]).
In the course of its preparation  we had useful discussions
on different related
topics with D. Akhiezer, F. Bogomolov, M. Brion,  R.O. Buchweitz,
F. Catanese, H. Kraft,
F. Kutzschebauch, S. Orevkov and V. Sergiescu. Their advice,
references and information
were very helpful. We are grateful to all of them.
The first named author would like to thank the Institut Fourier in
Grenoble and the second named author would like to thank the SFB 170
`Geometry and Analysis' in G\"{o}ttingen for their hospitality.

\section {Preliminaries}

\noin {\bf a) Background on plane algebraic curves}

\vs

\noin One says that a reduced curve $C$ in $\pr^2$ has {\it classical
singularities}
if all its singular points are nodes and ordinary cusps. It is called
{\it a Pl\"ucker curve} if both $C$ and the dual curve $C^*$
have only classical
singularities and  no flecnode, i.e. no flex at a
node\footnote{observe that the Pl\"{u}cker formulas are still valid
if the latter condition is omited, but in this case
 one must count separately the flexes and nodes which
 are coming from flecnodes or biflecnodes.}. We say
that
$C$ is {\it an immersed curve} if the normalization mapping
$\nu : C_{norm} \to C \hookrightarrow \pr^2$ is an immersion, or, which
is equivalent, if all the irreducible local analytic branches of $C$ are
smooth
(in particular, this is so if $C$ has only ordinary singularities
\footnote{i.e. singularities where all the local branches are smooth
and pairwise transversal.}).

 Let $C \s \pr^2$ be an irreducible curve of degree $d\ge 2$ and of
 geometric
 genus $g$. Then $d^* ={\rm deg}\,
 C^*$ (i.e. the class of $C$) is defined by the class formula (see
 Namba [25], (1.5.4))
\be d^* = 2(d+g-1)
- \sum\limits_{p \in {\rm sing}\, C} (m_p - r_p)\,, \ee
where $m_p = {\rm mult}_p C$ and $r_p$ is the number of irreducible analytic
branches of $C$ at $p$. Thus, $d^* \ge  2(d+g-1)$, where the equality holds
iff $C$ is an immersed curve.

We will need the following corollary of the genus formula (see Namba [25], (2.1.10)):
$$ 2g \leq  (d-1)(d-2) - \sum\limits_{p \in {\rm sing}\, C} m_p(m_p-1)$$
and $2g=(d-1)(d-2) -2 \delta$ for a nodal curve with $\delta$ nodes.
For reader's convenience we recall also the usual Pl\"ucker formulas:
$$g = 1/2 (d - 1)(d - 2) - \delta - \kappa = 1/2 (d^* - 1)(d^* - 2) -
b - f$$  $$d^* = d(d - 1) - 2\delta - 3\kappa\,\,\,\,\,{\rm and}\,\,\,\,\,d =
d^* (d^* - 1) - 2b - 3f$$ for a Pl\"ucker curve $C$ with $\delta$ nodes,
$\kappa$ cusps, $b$ bitangent lines and $f$ flexes. \vs

 Let $C \s \pr^2$ be an irreducible curve of degree $d \ge 2$ and let
 $\nu : C^*_{norm} \to C^*$ be the normalization of the dual curve.
 Following Zariski [38], p.307, p.326 and M. Green [16] (see also Dolgachev and Libgober [11]),
 consider the mapping $\rho_C \,: \,\pr^2 \to S^n C^*_{norm}$ of $\pr^2$
 into the $n$-th symmetric power of $C^*_{norm}$, where $n = {\rm deg}\,
 C^*$ and where $\rho_C (z) = \nu^* (l_z) \s  S^n C^*_{norm}$ (here $z \in
 \pr^2$ and $l_z \s \pr^{2*}$ is the dual line). It is easy to check that
 $\rho_C: \pr^2 \to S^n C^*_{norm}$
is a holomorphic embedding, which we call in the sequel {\it the Zariski
embedding}. We denote by $\pr_C^2$ the image $\rho_C (\pr^2)$ in $S^n
C^*_{norm}$, by $D_n$ the union of the diagonal divisors in
$(C^*_{norm})^n$ and by $\Delta_n = s_n (D_n) \s S^n C^*_{norm}$ the
discriminant divisor,
i.e. the ramification locus of the branched covering $s_n: (C^*_{norm})^n
\to  S^n C^*_{norm}$. Thus, we have the diagram
$$
\begin{picture}(1000,60)
\thicklines
\put(125,5){$C \subset \pr^2 $}
\put(180,20){$\rho_C$}
\put(180,5){$\hookrightarrow$}
\put(210,5){$\pr_C^2 \subset $}
\put(245,5){$S^n (C^*_{norm}) \supset \Delta_n $}
\put(250,45){${(C^*_{norm}})^n \supset D_n$}
\put(275,38){\vector(0,-1){20}}
\put(319,38){\vector(0,-1){20}}
\put(255,25){$s_n$}
\put(400,25){(2)}

\end{picture}
$$
It is easily seen that $C \subset \rho_C^{-1} (\Delta_n)$. Besides $C$, this
preimage may also contain some lines which we call {\it artifacts}.

To be more precise, denote by $L_C$ the union of the dual lines in
$\pr^2$ of the cusps of $C^*$ (by {\it a cusp} we mean an irreducible singular
local branch). Clearly, $L_C$ consists of the inflexional tangents of $C$ and
the cuspidal tangents at those cusps of $C$ which are not simple, i.e. which
can not be resolved by just one blow-up. Due to an analogy in tomography,
we call $L_C$ {\it the artifacts} of $C$. These artifacts arise
naturally as soon as $C^*$ is not immersed, namely we have
$$\rho_C^{-1}(\pr_C^2 \cap \Delta_n) = C \cup L_C \,\,.$$ Indeed,
a point $z \in \pr^2 \setminus C$ is contained in $ \rho_C^{-1} (\Delta_n)$
iff its dual line $l_z$ passes through a cusp of $C^*$.

 If $C \s \pr^2$ is a rational curve of degree $d > 1$, then
$C^*_{norm} \cong \pr^1,\,\,  S^n \pr^1  \cong \pr^n $, and hence the Zariski
embedding  $\rho_C$ embeds $\pr^2$ into $\pr^n \cong S^n \pr^1$,
where $n = {\rm deg}\,C^*$. The normalization map $ \nu \,:\,\pr^1
\to C^* \s \pr^2$ can be given as  $\nu = (g_0 : g_1 : g_2)\,$,
where $g_i (z_0 , z_1) = \sum\limits_{j=0}^n b^{(i)}_j z_0^{n-j}
z_1^j \,\,,\,\,i=0,1,2, $ are homogeneous polynomials of degree $n$ without
common factor.

If $x = (x_0 : x_1 : x_2) \in \pr^2$ and $l_x \s \pr^{2*}$ is the dual line,
then $\rho_C (x) = \nu^* (l_x) \in S^n \pr^1 = \pr^n$ is defined by the
equation $\sum\limits_{i=0}^2 x_i g_i (z_0 : z_1 ) =0$. Thus, $\rho_C (x)
= (a_0 (x) : \dots : a_n (x) )$, where $a_j (x) = \sum\limits_{i=0}^2 x_i
b^{(i)}_j$.

Therefore, in the case of a plane rational curve $C$ the Zariski embedding
$\rho_C \,:\,\pr^2 \to \pr^n$ is a linear embedding given by the
$3 \times (n+1)$--matrix $B_C := (b^{(i)}_j )$, $i=0,1,2,\,j=0,\dots,n$, and
its image $\pr^2_C = \rho_C (\pr^2)$ is a plane in $\pr^n$.

\vs

\noin {\bf b) On the Vieta map and the PGL$(2,\, \cz)$--action on
$\pr^n$}

\vs

\noin The symmetric power $S^n \pr^1$ is naturally identified with $\pr^n$ in
such a way that the canonical projection $s_n \,:\, (\pr^1 )^n \to S^n \pr^1$
coincides with {\it the Vieta ramified covering} given by
$$((u_1 : v_1),\dots,(u_n : v_n)) \longmapsto $$
$$\longmapsto (\prod\limits_{i=1}^n v_i)\,(1 : \sigma_1 (u_1 /v_1 ,
\dots, u_n /v_n )\,:\,\dots \,:\,\sigma_n (u_1 /v_1 ,\dots, u_n /v_n
))\,\,\,,$$
where $\sigma_i (x_1 ,\dots, x_n)\,,\,\,i=1,\dots,n$, are the elementary
symmetric polynomials. This is a Galois covering with the Galois group being
the n-th symmetric group $S_n$. With $z_i := (u_i : v_i) \in \pr^1 ,\,\,i=1,
\dots,n$, we have $s_n (z_1 ,\dots, z_n) = (a_0 : \dots : a_n)$, where
$z_i,\,i= 1,\dots, n,\,$ are the roots of the binary form $\,
\sum\limits_{i=0}^n a_i u^{n - i} v^i$ of degree $n$; see Zariski [38], p.252.

\vs

Note that the Vieta map $s_n\,:\,(\pr^1)^n \to S^n\pr^1 = \pr^n$
is equivariant with respect to the induced actions of the group
PGL$(2,\, \cz) = {\rm Aut}\,\pr^1$ on $ (\pr^1)^n$ and on $\pr^n$,
respectively. The branching divisors $D_n \subset (\pr^1)^n$
(the union of the diagonals)
resp. $\Delta_n \subset \pr^n$ (the discriminant divisor), as well as their
complements are
invariant under these actions. It is easily seen that for
$n \ge 3$ the orbit space of the PGL$(2,\, \cz)$--action on
$\pr^n \se \Delta_n$ is naturally isomorphic to the moduli space $M_{0,\,n}$
of the Riemann sphere with $n$ punctures. Denote by $\tilde {M}_{0,\,n}$
the quotient $((\pr^1)^n \se D_n)\, / \pr GL(2,\, \cz)$. We have the
following commutative diagram of equivariant morphisms

$$
\begin{picture}(800,70)
\unitlength0.2em
\thicklines
\put(54,25){$(\pr^1)^n \se D_n$}
\put(110,25){$\tilde M_{0,\,n}$}
\put(87,27){$\vector(1,0){15}$}
\put(60,5){$\pr^n \se \Delta_n$}
\put(87,6){$\vector(1,0){15}$}
\put(110,5){$M_{0,\,n}$}
\put(69,22){$\vector(0,-1){11}$}
\put(114,22){$\vector(0,-1){11}$}
\put(90,10){$\pi_n$}
\put(90,30){$\tilde \pi_n$}
\put(60,16){$s_n$}
\put(170,15){(3)}
\end{picture}
$$
The cross--ratios $\delta_i (z) = (z_1,\,z_2;\,z_3 ,\,z_i)$, where
$z = (z_1 ,\dots, z_n ) \in (\pr^1)^n$ and $4 \le i \le n$, define a
morphism $$\delta^{(n)} = (\delta_4 ,\,\dots, \,\delta_n )\,:\,(\pr^1)^n
\se D_n \to (\cz^{**})^{n - 3} \se D_{n - 3}\,,$$
where $\,\cz^{**} := \pr^1 \se \{0,\,1,\,\infty\}$. By the invariance of
cross--ratio $\delta^{(n)}$ is constant along the orbits of the action
of PGL$(2,\, \cz)$ on $(\pr^1 )^n \se D_n$. Therefore, it factorizes
through a mapping of the orbit space $\tilde M_{0,\,n} \to (\cz^{**})^{n - 3}
\se D_{n - 3}$. On the other hand, for each point $z \in (\pr^1)^n \se D_n$
its PGL$(2,\, \cz)$--orbit $O_z$ contains the unique point $z'$ of the
form $z' = (0,\,1,\,\infty,\,z'_4 ,\,\dots,\,z'_n )$. This defines a regular
section $\tilde M_{0,\,n} \to (\pr^1)^n \se D_n$, and its image coincides
with the image of the biregular embedding $$(\cz^{**})^{n - 3} \se D_{n - 3}
\ni u = (u_4 ,\,\dots ,\,u_n ) \longmapsto (0,\,1,\,\infty,\,u_4 ,\,\dots,
\,u_n ) \in (\pr^1)^n \se D_n \,\,.$$
This shows that the above mapping $\tilde M_{0,\,n} \to (\cz^{**})^{n - 3}
\se D_{n - 3}$ is an isomorphism. \vs

In the sequel we treat $\pr^n$ as the projectivized space of the binary
forms of degree $n$ in $u$ and $v$. For instance, $e_k = (0:\dots :0:1_k :0:
\dots :0) \in \pr^n$ corresponds to the forms $cu^{n - k}v^k$, where $c \in
\cz^*$. Denote by $O_q$ the PGL$(2,\, \cz)$--orbit of a point
$q \in \pr^n$; it is a smooth quasiprojective variety. If the form $q$ has
the roots $z_1,\,z_2 ,\dots$ of multiplicities $m_1,\, m_2 ,\, \dots$,
then we say that $O_q$ is an orbit of type $O_{m_1,\,m_2,\dots}$;
furthermore, even in the case when $O_q$ is not the only orbit of this type,
without abuse of notation we often write  $O_{m_1,\,m_2,\dots}$ for the
orbit  $O_q$ itself. Clearly, $O_{e_i} = O_{e_{n - i}},\,i=0,\dots, n$;
$O_{e_0} = O_n$ is the only one--dimensional orbit and, at the same time,
the only closed orbit; $O_{e_i} = O_{n-i,\,i} ,\,i=1,\dots,[n/2]$, are the
only two-dimensional orbits. Any other orbit $O_q = O_{m_1 ,\,m_2 ,\,m_3 ,
\dots}$ has dimension $3$ and its closure $\bar {O_q}$ is the union of the
orbits
$O_q, O_n$ and $O_{m_i,\,n-m_i} ,\,i=1,2,\dots$, which follows from
Aluffi and Faber [1], Proposition 2.1.
Furthermore, for any point $q \in \pr^n \se \Delta_n$, i.e. for any binary form
$q$ without multiple roots, its orbit $O_q = O_{1,\,1,\dots,1}$ is closed
in $\pr^n \se \Delta_n$, and its closure in $\pr^n$ is $\bar {O_q}
= O_q \cup S_1$, where $S_1 := O_n \cup O_{n-1,\,1} = \bar {O_q} \cap
\Delta_n$.
Therefore, any Zariski closed subvariety $Z$ of $\pr^n$ such that
${\rm dim}\,(O_q \cap Z) > 0$ must meet the surface $S_1$. These
observations yield the following lemma. \footnote {We are grateful to H.
Kraft for pointing out the approach used in the proof, and to M.
Brion for mentioning the paper Aluffi and Faber [1].}

\vs

\noin {\bf 2.1. Lemma.} {\it If a linear subspace $L$ in $\pr^n$ does not
meet the surface $S_1 = {\bar O}_{n-1,\,1} \s \Delta_n$, then it has at most
finite intersection with any of the orbits $O_q$, where  $q \in \pr^n \se
\Delta_n$. In particular, this is so for a generic linear subspace $L$ in
$\pr^n$
of codimension at least $3$.}

\vs

For instance, for a given  $k$--tuple of distinct points
$z_1,\dots, z_k \in \cz$, where $3\le k\le n$, consider the projective
subspace
$H_k = H_k(z_1,\dots, z_k) \subset \pr^n$, consisting of the
binary forms
of degree $n$ which vanish at $z_1,\dots, z_k$. Then, clearly, $H_k $
satisfies the above condition, i.e. it does not meet $S_1$.

\vs

\noin {\bf c) Background in hyperbolic complex analysis} \vs

\noin The next statement follows from Zaidenberg [35], Thms. 1.3, 2.5. \vs

\noin {\bf 2.2. Lemma.} {\it Let $C \s \pr^2$ be a curve such that the
Riemann surface ${\rm reg}\, C := C \se {\rm sing}\,C$ is hyperbolic
and $\pr^2 \se C$ is Brody hyperbolic, i.e. it does not contain any entire
curve. Then $\pr^2 \se C$ is Kobayashi complete hyperbolic and hyperbolically
embedded into $\pr^2$. The condition "${\rm reg}\, C$ is hyperbolic" is
necessary for
$\pr^2 \se C$ being  hyperbolically embedded into $\pr^2$. }

\vs

 We say that a complex space $X$ is {\it almost} resp. {\it weakly
 Carath\'eodory hyperbolic} if for any point $p \in X$ there exist only
 finitely many resp. countably many points $q \in X$
which cannot be separated from $p$ by bounded holomorphic functions.
It will be called {\it almost} resp. {\it weakly C--hyperbolic} if $X$ has
a  covering $Y \to X$, where $Y$ is almost resp. weakly Carath\'eodory
hyperbolic. Note that the universal covering $\tilde X$ of a  C--hyperbolic
complex manifold $X$ need not to be Carath\'eodory hyperbolic \footnote{F.
Kutzschebauch has constructed a corresponding example of a non--Stein domain
$X \s \s \cz^2$ (letter to the authors from 6.7.1995)}. At the same time,
it is weakly Carath\'eodory hyperbolic.

The next lemma is evident.

\vs

\noin {\bf 2.3. Lemma.}  {\it Let $f: Y \to X$ be a holomorphic  mapping of
complex spaces.
If $f$ is injective (resp. has finite resp. at most countable fibres) and
$X$ is  C--hyperbolic (resp. almost resp. weakly C--hyperbolic),
 then $Y$ is C--hyperbolic (resp. almost resp. weakly C--hyperbolic).}

\section{Proof of Theorem 1.1 and a generalization}

The next theorem gives sufficient conditions for the complement of an
irreducible plain curve $C$ and its artifacts $L_C$ to be C-hyperbolic.
In the particular
case when the dual curve $C^*$ is immersed (i.e. when $L_C = \emptyset$)
this leads to Theorem 1.1 of the Introduction.

\vs

\noin {\bf 3.1. Theorem.} {\it Let $C \s \pr^2$ be an irreducible curve of
genus $g$. Put $n = {\rm deg}\,C^*$ and $X = \pr^2 \se (C \cup L_C )$.

\noin a) If $g \ge 1$, then $X$ is C--hyperbolic.

\noin b) If $g = 0$, then $X$ is
almost C--hyperbolic if at least one of the following conditions is fulfilled:

\noin $b'$) $\,i(T_{p^*} A^* , \,A^* ;\,p^*) \le n - 2$ for any local analytic
branch $(A^* , \,p^* )$ of $C^*$;

\noin $b''$) $C^*$ has a cusp and it is not projectively equivalent to one of
the curves $(1 : g(t) : t^n ),\,\,(t : g(t) : t^n )$, where $g \in \cz [t]$
and $ {\rm deg}\,g \le n - 2$.

\noin c) Under any of the assumptions of ($a$), ($b'$), ($b''$) $X$ is complete
hyperbolic and hyperbolically embedded into $\pr^2$ iff the curve
reg$\,(C \cup L_C) = (C \cup L_C) \setminus $ sing $(C \cup L_C)$ is
hyperbolic. }

\vs

\noin Here $i(.,.;.)$ stands for the local intersection multiplicity.

\vs

The last statement ($c$) easily follows from the previous ones in view of
Lemma 2.2 and the subsequent remark. Before passing to the proof of ($a$) and
($b$) let us make the following observations.

\vs

\noin {\it Remark.} Observe that under the conditions of Theorem 1.1 the curve
reg$\,C$ is hyperbolic. Indeed, the dual of an immersed curve can not be
smooth; therefore, this is true as
soon as $g \ge 1$, i.e. under the condition in a). This is also true if $C^*$
is a generic rational curve of degree $n \ge 5$, as it was supposed in
(b). More generally, let $C$ be a rational curve of degree $d $ such that
the dual  $C^*$ is an immersed curve of degree $n > 2$. Then by class formula
(1) we have:
$d = 2(n-1)$ and
$$ \sum\limits_{p \in {\rm sing}\, C} (m_p - r_p) = 2(d-1)-n = 3(n-2)
\ge 3\,.$$ Thus, $C$ has at least three cusps and therefore, reg$\,C$ is
hyperbolic.

Notice also that condition ($b'$) is fulfilled for the dual of a generic
 rational curve of degree $n \ge 5$. This ensures that, indeed,
Theorem 1.1 follows from Theorem 3.1.

\vs

\noin {\bf Proof of Theorem 3.1.$a$}

\vs

\noin Let $\rho_C\,:\,\pr^2 \to S^n (C^*_{norm})$ be the Zariski embedding
introduced in section 2.a). The covering $s_n: (C^*_{norm})^n \se D_n \r
S^n(C^*_{norm}) \se \Delta_n$ is non--ramified. Thus, we have the commutative
diagram
$$
\begin{picture}(1000,60)
\thicklines
\put(127,5){$X$}
\put(130,45){$Y$}
\put(132,38){\vector(0,-1){20}}
\put(105,25){${\tilde s}_n$}

\put(175,58){${\tilde \rho}_C$}
\put(175,20){$\rho_C$}
\put(175,45){$\hookrightarrow$}
\put(175,5){$\hookrightarrow$}

\put(220,5){$S^n(C^*_{norm} ) \se \Delta_n \,\,\,\,$}
\put(225,45){$(C^*_{norm})^n \se D_n$}
\put(260,38){\vector(0,-1){20}}
\put(230,25){$s_n$}
\put(420,25){(4)}
\put(310,45){$\hookrightarrow$}
\put(340,45){${(C^*_{norm})}^n$}

\end{picture}
$$
where ${\tilde s}_n \,:\,Y \to X$ is the induced covering. If
genus $g(C^*) \ge 2$, then $(C^*_{norm} )^n$ has the polydisc
$\D ^n$ as the universal covering. Passing to the induced covering
$Z \to Y$ we can extend (4) to the diagram

$$
\begin{picture}(1000,80)
\thicklines
\put(155,75){$Z$}
\put(158,69){\vector(0,-1){18}}
\put(155,40){$Y$}
\put(158,36){\vector(0,-1){18}}
\put(153,5){$X$}
\put(135,25){${\tilde s}_n$}
\put(200,75){$\hookrightarrow$}
\put(200,40){$\hookrightarrow$}
\put(200,5){$\hookrightarrow$}
\put(200,55){${\tilde \rho}_C$}
\put(200,20){$\rho_C$}
\put(253,75){$U^n$}
\put(258,69){\vector(0,-1){18}}
\put(246,40){$({C^*_{norm})}^n$}
\put(258,36){\vector(0,-1){18}}
\put(247,5){$S^n (C^*_{norm})$}
\put(264,25){$s_n$}
\put(420,40){(5)}
\end{picture}
$$
Being a submanifold of the polydisc $Z$ is Carath\'eodory hyperbolic, and
so $X$ is C--hyperbolic. Therefore, we have proved Theorem 3.1.$a$ in the
case $g \geq 2$.\vs

Next we consider the case $g=1$. Denote $E = C^*_{norm}$. Note that both
$E^n \se D_n$ and $S^nE \se \Delta_n$ are not C--hyperbolic or even hyperbolic,
and so we can not apply the same arguments as above.

Represent $E$ as $E = J(E) = \cz / \Lambda_{\omega}$, where
$\Lambda_{\omega}$ is the lattice generated by $1$ and $\omega
\in \cz_{+}$ (here $ \cz_{+} := \{z \in \cz\,|\,{\rm Im}z > 0\}$). By
Abel's Theorem we may assume this identification of
$E$ with its jacobian $J(E)$ being chosen in such a way that the
image $\rho_C (\pr^2)$ is contained in the hypersurface $s_n (H_0) =
\phi_n^{-1} ({\bar 0}) \cong \pr^{n-1} \s S^n E$, where $$H_0 :=
\{z=(z_1,\dots, z_n) \in E^n\,|\,\sum\limits_{i=1}^n z_i = 0\}$$ is an
abelian subvariety in $E^n$ and $\phi_n\,:\,S^n E \to J(E)$ denotes the
n-th Abel--Jacobi map. The universal covering $\tilde{H}_0$ of $H_0$
can be identified with the hyperplane $\sum\limits_{i=1}^n x_i = 0$
in $\cz^n = \tilde{E}^n$.

Consider the countable families $\tilde{D}_{ij}$ of parallel affine
hyperplanes in $\cz^n$ defined by the conditions $x_i - x_j \in
\Lambda_{\omega}\,,\,\,\,i,j = 1,\dots,n,\,\,i<j$.

\vs

\noin {\it Claim. The domain $\tilde{H}_0 \se \bigcup\limits_{i=1}^{n-1}
\tilde{D}_{i,i+1}$ is biholomorphic to $(\cz \se \Lambda_{\omega})^{n-1}$.}

\vs

\noin Indeed, put $y_k := (x_k - x_{k+1})\,|\,\tilde{H}_0\,,\,\,i =
1,\dots,n-1$.
It is easily seen that $(y_1,\dots,y_{n-1})\,:\, \tilde{H}_0 \to \cz^{n-1}$
is a linear isomorphism whose restriction yields a biholomorphism as in the
claim.

The universal covering of  $(\cz \se \Lambda_{\omega})^{n-1}$ is the
polydisc $\D^n$, and so $(\cz \se \Lambda_{\omega})^{n-1}$  is C--hyperbolic.
Put $\tilde{D}_n:= \bigcup\limits_{i,j=1,\dots,n} \tilde{D}_{ij}$. The open
subset $\tilde{H}_0 \se \tilde{D}_n$ of $\tilde{H}_0 \se
\bigcup\limits_{i=1}^{n-1} \tilde{D}_{i,i+1} \cong (\cz \se
\Lambda_{\omega})^{n-1} $ is also C--hyperbolic (see (2.3)).

Denote by $p$ the universal covering map $\cz^n \to (\cz /\Lambda_{\omega})^n$.
The restriction $$p\,|\,\tilde{H}_0 \se \tilde{D}_n\,:\,\tilde{H}_0 \se
\tilde{D}_n  \to H_0 \se D_n \s E^n \se D_n$$ is also a covering map.
Therefore, $H_0 \se D_n$ is C--hyperbolic, and so $s_n (H_0) \se \Delta_n$ is
C--hyperbolic, too. Since  $\rho_C\,|\, X \,:\, X \to s_n (H_0) \se \Delta_n$
is a holomorphic embedding, by Lemma 2.3 $X$ is C--hyperbolic. \qed

\vs

\vs

\noin {\bf Proof of Theorem 3.1.$b'$}

\vs

It consists of the next two lemmas. We freely use the notation from
sect. 2.a, 2.b. Remind, in particular, that for a rational curve
$C \subset \pr^2$ with
the dual $C^*$ of degree $n$, the plane
$\pr^2_C = \rho_C(\pr^2) \subset \pr^n$
is the image of $\pr^2$ under the Zariski embedding. The surface $S_1 \subset
\Delta_n
\subset \pr^n$ is the orbit closure ${\bar O}_{n-1,\,1}$.

\vs

\noin {\bf 3.2. Lemma.}  {\it The complement $X = \pr^2 \se (C \cup L_C)$,
where $C \subset \pr^2$ is a rational curve, is almost C--hyperbolic
whenever $\pr^2_C \cap S_1 = \emptyset$.}

\vs

\noin {\it Proof.} Consider the following commutative diagram of morphisms:

\begin{picture}(250,80)
\unitlength0.2em
\thicklines
\put(33,5) {$X$}
\put(35,25){$Y$}
\put(36,22){\vector(0,-1){11}}
\put(27,16){${\tilde s}_n$}
\put(44,27){{\vector(1,0){14}}}
\put(49,30){${\tilde\rho}_C$}
\put(44,6){\vector(1,0){14}}
\put(49,10){$\rho_C$}
\put(66,5){$\pr^n \se \Delta_n$}
\put(62,25){$(\pr^1)^n \se D_n$}
\put(76,22){\vector(0,-1){11}}
\put(128,22){\vector(0,-1){11}}
\put(124,5){$M_{0,\,n}$}
\put(67,16){$s_n$}
\put(170,15){(6)}
\put(92,27){{\vector(1,0){14}}}
\put(94,6){{\vector(1,0){20}}}
\put(95,30){${\tilde \pi}_n$}
\put(100,10){$\pi_n$}
\put(110,25){$ (\cz^{**})^{n - 3} \se D_{n - 3} \hookrightarrow
(\cz^{**})^{n - 3}$}
\end{picture}

\noin where ${\tilde s}_n \,:\,Y \to X$ is the induced covering (cf. (4)).
{}From Lemma 2.1 it follows that the mapping $\pi_n \circ \rho _C \,:\,X \to
M_{0,\,n}$ has finite fibres. Hence, the same is valid for the mapping ${\tilde
\pi}_n \circ {\tilde\rho}_C \,:\, Y \to (\cz^{**})^{n - 3} \se D_{n - 3}$. By
Lemma 2.3 $Y$, and thus also $X$, are almost C--hyperbolic. \qed

\vs

\noin {\bf 3.3. Lemma.} {\it Let $C \s \pr^2$ be a rational curve.
Put $n = {\rm deg}\,C^*$. Then $\pr^2_C \cap S_1 = \emptyset$ iff the condition
($b'$) is fulfilled, i.e. iff $\,i(T_{p^*} A^* , \,A^* ;\,p^*) < n - 1$ for any
local analytic branch $(A^*, p^* )$ of the dual curve $C^*$.}

\vs

\noin{\it Proof.} By definition of the Zariski embedding,
$q \in \pr^2_C \cap S_1$ iff, after passing to normalization $\nu :\, \pr^1 \to
C^*$ and identifying $\pr^2$ with its image $\pr^2_C$ under $\rho_C$, the dual
line $l_q \s \pr^{2*}$ cuts out on $C^*$ a divisor of the form $(n - 1)a + b$,
where $a, b \in \pr^1$. Then $p^* := \nu (a) \in C^*$ is the center of a local
branch $A^*$ of $C^*$ which violates the condition in ($b'$). The converse is
evidently true. \qed

\vs

\noin {\it Remarks.} 1. If the dual curve $C^*$ has only ordinary cusps and
flexes and $n = {\rm deg}\,C^* \ge 5$, then $\pr^2_C \cap S_1 = \emptyset$.
Indeed, in this case $i(T_{p^*} A^* , \,A^* ;\,p^*) \le 3 < n - 1$ for any
local analytic branch $(A^*, p^* )$ of $C^*$, and so the result follows from
Lemma 3.3.

\vs

\noin 2.  If $C^*$ has a cusp $(A^*, p^* )$ of multiplicity $n - 1$, then
$\rho_C (l_{p^*}) \s \pr^2_C \cap S_1$, where $l_{p^*} \s L_C \s \pr^2$ is
the dual line of the point $p^* \in \pr^{2*}$. Indeed, for any point
$q \in l_{p^*}$ its dual line $l_q \s \pr^{2*}$ passes through $p^*$,
and hence we have, as above, $\rho_C (q) \in \pr^2_C \cap S_1$.

\vs

Next we give an example where both of the conditions $b',\,b''$ of
Theorem 3.1 are violated.

\vs

\noin {\bf 3.4. Example.} Let $C^* = (p(t) : q(t) : 1)$ be a parametrized plane
rational curve, where $p, \,q \in \cz [t]$ are generic polynomials of degree
$n$ and $n - 1$, respectively. Thus, $C^*$ is a nodal curve of degree $n$
which is the projective closure of an affine plane polynomial curve with one
place at infinity at the point $(1 : 0 : 0)$ which is a smooth point of $C^*$.
The line at infinity $l_2 = \{x_2 = 0 \}$ is an inflexional tangent of $C^*$
(of order $n - 2$). By Lemma 3.3, $\pr^2_C \cap S \neq \emptyset$.
Therefore, Lemma 3.2 is not applicable. We do not know whether the complement
$\pr^2 \se C$
of the dual $C$ of $C^*$ in this example is C--hyperbolic or not.

\section{Projective duality and $\cz^*$--actions}

The proof of Theorem 3.1.$b''$ is based on a different idea. It needs
certain preparations, which is the subject of this section; the proof
is done in the next one.

\vs

\vs

\noin {\bf a) Veronese projection, Zariski embedding and projective duality}

\vs

Let $C \subset \pr^2$ be a rational curve
with the dual $C^*$ of degree $n$, and
$\rho_C\,:\,\pr^2 \hookrightarrow \pr^2_C \subset \pr^n$ be the Zariski
embedding.
The dual map ${\rho_C}^* \,:\,\pr^{n*} \to \pr^{2*}$ given by the transposed
matrix $\,^t B_C$ (see sect. 2.a) defines a linear projection with center
$N_C := {\rm Ker\,} \,^t B_C \s \pr^{n*}$ of codimension $3$. The curve
 $C^*$ is the image under this projection of the rational normal curve
$C_n^* = (z_0^n : z_0^{n-1} z_1 : \dots : z_1^n ) \s \pr^{n*}$ (see Veronese [33], p.208),
i.e. $\rho_C^* (C_n^*) = C^* $.
  Furthermore, $C_n^*$  is the image of $\pr^1 \cong C^*_{norm}$ under the
embedding $i \,:\, \pr^1 \hookrightarrow \pr^{n*}$ defined by the complete
linear system $|H| = |n(\infty)| \cong \pr^n$. The composition
 $\nu = {\rho_C}^* \circ i \,:\, \pr^1 \to C^* \s \pr^{2*}$ is the
normalization map.

The rational normal curve $C_n^* \s \pr^{n*}$ and the discriminant hypersurface
$\Delta_n \s \pr^{n}$ are dual to each other. This yields the following
duality:
\begin{center}
\begin{picture}(500,70)
\put(150,55){$(\pr^2 ,\,C \cup L_C)\,\,\,\,\,\,\, \hookrightarrow
\,\,\,\,\,\,\, (\pr^n ,\,\Delta_n )$}
\put(155,5){$(\pr^{2*} ,\,C^* ) \,\,\,\,\,\,\,\,\,\, \longleftarrow
\,\,\,\,\,\,\, (\pr^{n*} ,\,C_n^* )$}
\put(175,30){$\updownarrow$}
\put(291,30){$\updownarrow$}
\put(238,20){$\rho^*_C$}
\put(240,68){$\rho_C$}
\end{picture}
\end{center}
\noin To describe this duality in more details, fix a point
$q=(z_0^n:z_0^{n-1}z_1:...:z_1^n) \in C_n^* \s \pr^{n*}$, and let
$$F_q C^*_n = \{ T^0_q C_n^* \s T^1_q C_n^* \s \dots \s T^{n - 1}_q
C_n^* \s \pr^{n*} \}$$ be the flag of osculating subspaces to $C_n^*$ at $q$,
where ${\rm dim}\,T^k_q C_n^* = k,\,\,T^0_q C_n^* = \{q\}$ and $T^1_q C_n^*
= T_q C_n^*$ is the tangent line to $C_n^*$ at $q$ (see Namba [25], p.110). For
instance, for $q = q_0 = (1 : 0 : \dots : 0) \in C_n^*$ we have $T^k_q C_n^*
= \{x_{k + 1} = \dots = x_n = 0 \} \s \pr^{n*}$.

The dual curve $C_n \s \pr^n$ of $C_n^*$ is in turn projectively equivalent
to a rational normal curve; namely, $$C_n = \{ p \in \pr^n \,|\, p = q^* =
(z_1^n : -nz_0 z_1^{n - 1} : \dots : (-1)^k  { n \choose k } z_0^k z_1^{n - k}
: \dots : (-1)^n z_0^n )\} = O_n\,.$$
Furthermore, the dual flag $F_q^{\perp} = \{\pr^n \supset H^{n - 1}_q \supset
\dots \supset H^0_q \}$, where $H^{n - k}_q := (T^{k - 1}_q C_n^*)^{\perp}$, is
 the flag of osculating subspaces
$F_p C_n = \{T^{k - 1}_p C_n \}_{k = 1}^n $ of the dual rational normal curve
$C_n \s \pr^n$. An easy way to see this is to observe that at the dual points
$q_0 = (1 : 0 : \dots : 0 ) \in C_n^*$ and $p_0 = q_0^* = (0 : \dots : 0 : 1 )
\in C_n$ both flags consist of coordinate subspaces, and then to use ${\rm
Aut}\,\pr^1$-homogeneity.

The points of the osculating subspace $H^k_q = T^k_p C_n$ correspond to the
binary forms of degree $n$ for which $(z_0 : z_1) \in \pr^1$ is a root of
multiplicity at least $n - k$. In particular, $H_q^{n - 2} = (T_q C^*_n
)^{\perp}$ consists of the binary forms which have $(z_0 : z_1)$ as a multiple
root. Therefore, the discriminant hypersurface $\Delta_n$ is the union of
these linear subspaces $H_q^{n - 2} \cong \pr^{n - 2}$ for all $q \in C_n^*$,
and thus it is the dual hypersurface of the rational normal curve $C_n^*$, i.e.
each of its points corresponds to a hyperplane in $\pr^{n*}$ which contains a
tangent line of $C_n^*$. At the same time, $\Delta_n$ is the developable
hypersurface of the $(n-2)$--osculating subspaces  $H_q^{n - 2} = T_p^{n - 2}
C_n$ of the dual rational normal curve $C_n \s \Delta_n$; here $T_p^{n - 2}
C_n \cap C_n = \{p\}$.

Let $d_{ij} \cong \pr^1$ be the diagonal of $\pr^1_i \times \pr^1_j$. The
decomposition $D_{ij} = d_{ij} \times (\pr^1 )^{n-2}$ of the diagonal
hyperplane $D_{ij} \s D_n$ may be regarded as the trivial fibre bundle
$D_{ij} \to \pr^1$ with the fibre $(\pr^1 )^{n-2}$. The subspaces $H_q^{n - 2}
\s \Delta_n$ are just the images of the fibres under the Vieta map
$s_n\, :\,(\pr^1)^n \to \pr^n$. Moreover, the restriction of $s_n$ to a fibre
yields the Vieta map $s_{n-2}\, :\,(\pr^1)^{n-2} \to \pr^{n-2}$. The dual
rational normal curve $C_n \s \pr^n$ is the image $s_n (d_n )$ of the diagonal
line $d_n := \bigcap_{i,\,j} D_{ij} = \{z_1 = \dots = z_n \} \s (\pr^1 )^n$.

\vs


By duality we have $N_C={\rm Ker}\,
{\rho_C}^* = ({\rm Im}\,\rho_C )^{\perp}$, i.e. $N_C = (\pr^2_C )^{\perp}$.
Therefore,
$$\pr^2_C = N_C^{\perp} = \bigcap\limits_{x^* \in N_C } {\rm Ker}\,x^* = \{ x
\in \pr^n \,|\, <x, x^* > = 0 \,\,{\rm for\,\, all\,\,}x^* \in N_C \} \,\,\,.$$

A point $q$ on the rational normal curve $C_n^* \s \pr^{n*}$ corresponds to a
cusp of $C^*$ under the projection ${\rho_C}^*$ iff the center $N_C$ of the
projection meets the tangent developable $TC_n^* = S_1$, which is a ruled
surface in $\pr^{n*}$, in some point $x_q^*$ of the tangent line $T_q C_n^*$
(see Piene [28]). In this case it meets $T_q C_n^*$ at the only point $x_q^*$,
because otherwise $N_C$ would contain $T_q C_n^*$ and thus also the point $q$,
which is impossible since ${\deg}\,C^* = {\deg}\,C_n^* = n$.

Let $B$ be a cusp (i.e. a singular local analytic branch) of $C^*$ centered at
the point $ q_0 = {\rho_C}^* (q)$. It
corresponds to a local branch of $C_n^*$ at the point $q \in C_n^*$ under the
normalizing projection ${\rho_C}^* \,:\,C_n^* \to C^*$. Define $L_{B, q_0} :=
{\rm Ker}\,x_q^* \s \pr^n$ to be the dual hyperplane of the point $x_q^*
\in N_C \cap T_q C_n^*$. Since $x_q^* \in N_C$, this  hyperplane $L_{B,  q_0}$
contains the image $\pr^2_C = \rho_C (\pr^2 )$. This yields a
correspondence between the cusps of $C^*$
and certain hyperplanes in $\pr^n$ containing the plane $\pr^2_C$. From the
definition it follows that $L_{B,q_0}$ contains also the dual linear space
$H_q^{n - 2} = (T_q C_n^*)^{\perp} \s \Delta_n$ of dimension $n - 2$.
Since the plane $\pr^2_C$ is not contained in $\Delta_n$, we have
$L_{B,q_0} = {\rm span}\, (\pr^2_C ,\,H_q^{n - 2} )$. It is easily
seen that the intersection $\pr^2_C \cap H_q^{n - 2}$ coincides with
the tangent line $l_{q_0} \s L_C$ of $C$, which is dual to the cusp $q_0$ of
$C^*$. Thus, the artifacts
$L_C$ of $C$ are the sections of $\pr^2_C$ by those osculating linear
subspaces $H_q^{n - 2} \s \Delta_n$ for which  $q$ is a cusp of $C^*$;
any other subspace $H_{q'}^{n - 2}$ meets the plane $\pr^2_C$ in one point
of $C$ only.

\vs

In what follows by {\it a special normalization} of the dual rational curve
$C^*$ we mean a normalization $\nu : \pr^1 \r C^* \s \pr^{2*}$ given in an
affine chart in $\pr^1$ as $\nu = (h_0(t) : h_1(t) : h_2(t)+t^n)$, where
$h_i \in \cz [t]$ and deg$\,h_i \leq n-2\,,\,\,i=0,1,2$. Such a curve
$C^*$ has a cusp $B$ at the point $q_0 = (0 : 0 : 1)$ which corresponds to $t
= \infty$. We will see below that $L_{B, q_0} = {\bar A}_1$, where
$${\bar A}_1 := \{(a_0 : \dots : a_n) \in \pr^n \,|\,a_1 = 0 \}\,\,.$$
Clearly, the preimage ${\bar H}_0 := s_n^{-1} ({\bar A}_1) \s (\pr^1)^n$ is
the closure of the linear hyperplane in $\cz^n$
$${H}_0 := \{z = (z_1,\dots, z_n ) \in \cz^n \,|\, \sum\limits_{i=1}^n
z_i =0 \} \,.\,\,$$
Note that the choice of normalization of $C^*$ is defined up to the PGL$(2,\,
\cz)$--action on $\pr^1$, and the induced
PGL$(2,\, \cz)$--action on $\pr^n$ affects the Zariski embedding. The
next lemma ensures the existence of special normalizations.

\vs

\noin {\bf 4.1. Lemma.} {\it Let $C^* \s \pr^{2*}$ be a rational curve
of degree $n$
with a cusp $B$ centered at the point $q_0  = (0:0:1) \in C^*$, and let
$L_{B, q_0} \s \pr^n$ be the corresponding hyperplane which contains the
plane $\pr^2_C = \rho_C (\pr^2 )$. Then $C^*$ admits a special normalization,
and under this normalization we have $L_{B, q_0} = {\bar A}_1$, where ${\bar
A}_1$ is as above.}

\vs

\noin {\it Proof.} The normalization
$\nu \,:\,\pr^1 \cong C^*_{norm} \to C^* \hookrightarrow \pr^2$ can be chosen
in such a way that the cusp $B$ corresponds to the local branch of $\pr^1$ at
$\infty= (1:0) \in \pr^1$, and so $\nu (\infty) = q_0$. If $\nu = (g_0 : g_1 :
g_2)\,$ is given by a triple of homogeneous polynomials
$g_i (z_0 , z_1) = \sum\limits_{j=0}^n b^{(i)}_j z_0^{n-j} z_1^j
\,\,,\,\,i=0,1,2, $ of degree $n$, then since $\nu (\infty) = q_0=(0 : 0 : 1)$
we have
${\rm deg}_{z_0}\,g_0 < n\,,\,{\rm deg}_{z_0}\,g_1 < n\,,\,{\rm deg}_{z_0}\,
g_2
 = n$, i.e. $b^{(0)}_0 =b^{(1)}_0 =0\,,\,b^{(2)}_0 \neq 0$. Performing the
Tschirnhausen transformation
$$\pr^1  \ni (z_0 \,:\,z_1) \longmapsto (z_0 - {b_1^{(2)} \over nb_0^{(2)}}
z_1 \,: \, z_1)    \in \pr^1$$
we may assume, furthermore, that $b_1^{(2)} =0$.

\vs

\noin {\it Claim 1. The normalization $\nu$ as above is a special one, and
the image
$\pr^2_C = \rho_C (\pr^2 )$ is contained in the hyperplane ${\bar A}_1$.}

\vs

Indeed, since $C^*$ has a
cusp at  $q_0$, we have $(g_0 / g_2)'_{z_1} = (g_1 / g_2)'_{z_1} = 0$ at
the point $(1:0) \in \pr^1$, i.e. $(g_0)'_{z_1} = (g_1)'_{z_1} = 0$ when $z_1
= 0$. This means that ${\rm deg}_{z_0}\,g_0 < n-1\,,\,{\rm deg}_{z_0}\,g_1 <
n-1$, i.e. $b^{(0)}_1 =b^{(1)}_1 =0$. And also $b^{(2)}_1 = 0$,
as it has been achieved above by making use of the Tschirnhausen
transformation.

Since $b_1^{(i)} =0\,,\,i=0,1,2$, we have $a_1 (x) \equiv 0$. Therefore,
$\rho_C (x) \in {\bar A}_1$ for any $x \in \pr^2$, which proves the claim.
\qed

\vs

\noin {\it Claim 2. The dual space $H_q^{n - 2}$ to $T_qC_n^*$ is contained
in ${\bar A}_1$.}

\vs

Indeed, since $\nu ( \infty) = q_0$ and $\nu = \rho_C^* \circ i$ with
$i: \pr^1 \r C_n^* \s \pr^{n*}$ we get $q=(1:0:...:0)$. Thus, by the
preceding considerations
the subspace $H_q^{n - 2} = (T_qC_n^*)^{\perp}$ is given by the equations
$\{a_0=a_1=0 \} $, and hence it is contained in ${\bar A}_1$. \qed

\ss

As before, we have $L_{B,q_0} = {\rm span}\, (\pr^2_C ,\,H_q^{n - 2} )$.
Therefore, from Claims 1 and 2 we obtain $L_{B,q_0} = {\bar A}_1$. \qed

\vs

\vs

\noin {\bf b) Monomial and quasi--monomial rational plane curves}

\vs

We will use the following terminology. By {\it a parametrized rational plane
curve} we mean a rational curve $C$ in $\pr^2$ with a fixed normalization
$\pr^1 \to C$ of it. {\it A parametrized monomial} resp. {\it a parametrized
quasi--monomial plane curve} is a parametrized rational plane curve such that
all resp. two of its coordinate functions are monomials; the image curve
itself is then called {\it monomial} resp.\  {\it quasi--monomial}.

Recall that if $C = (g_0 : g_1 : g_2 )$, where $g_i
\in \cz [t] ,\,\,i=0,1,2$, is a parametrized rational plane curve, then the
dual curve $C^*$ has (up to canceling the common factors)
the parametrization $C^* = (M_{12} : M_{02} : M_{01} )$, where $M_{ij}$
are the $2 \times 2$--minors of the matrix
$$ \left( \begin{array}{ccc}
g_0 & g_1 & g_2 \\
g'_0 & g'_1 & g'_2
\end{array} \right)$$
The equation of $C$ can be written as $\frac{1}{x_2^d }
 {\rm Res} \,(x_0 g_2 - x_2 g_0 , x_1 g_2 - x_2 g_1 ) = 0$, where $d = {\rm
deg}\,C$ and ${\rm Res}$ means resultant (see e.g. Aure [2], 3.2).

Note that a linear pencil of monomial curves $C_{\mu} = \{\alpha x_0^l +
\beta x_1^{l-k} x_2^k =0\}$, where $\mu = (\alpha : \beta) \in \pr^1$, is
self--dual, i.e. the dual curve of a monomial one is again monomial and
belongs to the same pencil. In contrast, the dual
curve to a quasi--monomial one is not necessarily projectively equivalent
to a quasi--monomial curve (recall that two plane curves $C,\,C'$ are {\it
projectively equivalent} if $C' = \alpha ( C)$ for some $\alpha \in \pr
{\rm GL} (3;\, \cz ) \cong {\rm Aut}\,\pr^2$). The simplest example is the
nodal cubic $C = \{(x_0 : x_1 : x_2 ) = (t : t^3 : t^2 -1 )\}$. Indeed, its
dual curve is a quartic with three cusps; but a quasi--monomial curve may
have at most two cusps.

Observe that, while the action of the projective group $PGL(3, \cz)$ on $\pr^2$
does not affect the image $\pr^2_C = \rho_C (\pr^2 ) \s \pr^n = S^n \pr^1$,
 the choice of the normalization $\pr^1 \to C^*$, defined up to the action of
the group $PGL(2, \cz ) = {\rm Aut} \pr^1$, usually does. This is why in the
next lemma we have to fix the normalization of a rational plane
curve $C$. This automatically fixes a normalization of its dual curve $C^*$,
and vice versa.

Clearly, projective equivalence between parametrized curves is a stronger
relation than just projective equivalence between underlying projective
curves themselves.

\vs

\noin {\bf 4.2. Lemma.} {\it A parametrized rational plane curve
$C^* \s \pr^{2 *}$ of degree $n$ is projectively equivalent to a parametrized
monomial resp. quasi--monomial curve iff $\pr^2_C \s \pr^n$ is a coordinate
plane resp. contains a coordinate axis. This axis is unique iff $C^*$ is
projectively equivalent to a parametrized quasi--monomial curve, but not to a
monomial one.}

\vs

\noin {\it Proof.} Let $\nu \,:\,t \longmapsto (at^k : bt^m : g(t))$, where
$a,\,b \in \cz^* ,\,\,g \in \cz [t]$ and $t=z_0 / z_1 \in \pr^1$, define a
parametrized quasi--monomial curve $C^* \s \pr^{2 *}$ of degree $n$. Denote
$e_k = (0 : \dots : 0 : 1_k : 0 : \dots : 0) \in \pr^n$. Then $\rho_C$ is
given by the matrix $B_C = (b^{(0)},\,b^{(1)},\,b^{(2)}) = (ae_{n-k} ,\,
be_{n-m} ,\,b^{(2)})$, and therefore $\pr^2_C = \rho_C (\pr^2 ) = {\rm
span}\,(b^{(0)},\,b^{(1)},\,b^{(2)})$ contains the coordinate axis
$l_{n-k,\,n-m}$, where $l_{i,j} := {\rm span}\,(e_i ,\,e_j ) \s \pr^n$.

If $C^*$ is a parametrized monomial curve, i.e. if $g(t) = ct^r$, where $c
\in \cz^*$, then clearly $\pr^2_C$ is the coordinate plane
$\pr_{n-k,\,n-m,\,n-r}
:=  {\rm span}\,(e_{n-k},\,e_{n-m},\,e_{n-r})$. Actually, up to a permutation
there should be $0 = r < m < k = n$ and gcd$(m, \,n) = 1$; thus, $\pr^2_C =
\pr_{0,\,n-m,\,n}$ is a rather special coordinate plane.

Since the projective equivalence of parametrized plane curves does not affect
the $\pr^2_C$, this yields the first statement of the lemma in one direction.

Vice versa, suppose that  $\pr^2_C$ coincides with the coordinate plane
$\pr_{n-k,\,n-m,\,n-r}$. Performing a suitable linear coordinate change in
$\pr^{2 *}$ we may assume that $b^{(0)} = e_{n-k} , \,b^{(1)} = e_{n-m}
,\,b^{(2)} = e_{n-r}$, i.e. that $\nu (t) = (t^k : t^m : t^r)$. Therefore,
in this case the parametrized curve $C^*$ is projectively equivalent to a
monomial curve.

Suppose now that $\pr^2_C$ contains the coordinate axis $l_{n-k ,\,n-m}$.
Performing as above a suitable linear coordinate change in $\pr^{2 *}$ we may
assume that $b^{(0)} = e_{n-k} , \,b^{(1)} = e_{n-m}$, and so $\nu (t) = (t^k
: t^m : g(t))$. In this case $C^*$ is projectively equivalent to a
parametrized quasi--monomial curve. This proves the first assertion of the
lemma.

Let $C^* = (at^{n-k} : bt^{n-m} : g(t))$ be a  parametrized quasi--monomial
curve which is not projectively equivalent to a monomial one. Then as above
$\pr^2_C \supset l_{k,\,m}$, and this is the only coordinate axis contained in
$\pr^2_C$ (indeed, otherwise $\pr^2_C$ would be a coordinate plane, that has
been excluded by our assumption). The opposite statement is evidently true.
This concludes the proof. \qed

\vs

\vs

\noin {\bf c) $\,\,\,\cz^*$--actions}

\vs

The natural $\cz^*$--action on $\pr^1$ induces (via the
Aut$\pr^1$--representations as in sect. 2.b) the following
$\,\cz^*$--actions on $(\pr^1)^n$ resp. on $\pr^n = S^n \pr^1$:
$${\tilde{G}} \,:\,\cz^* \t  (\pr^1)^n \ni (\lambda,\,((u_1 : v_1),\dots,
 (u_n : v_n))) \longmapsto ((\lambda u_1 : v_1), \dots, (\lambda u_n : v_n))
\in (\pr^1)^n$$
resp.
$$G \,:\,\cz^* \t \pr^n \ni (\lambda,\, (a_0 : a_1 : \dots : a_n ))
\longmapsto (a_0 : \lambda a_1 : \lambda^2 a_2 : \dots : \lambda^n a_n ) \in
\pr^n\,.$$
The Vieta map $s_n : (\pr^1)^n \to \pr^n$ (see sect. 2.b) is equivariant with
respect  to these $\,\cz^*$--actions and its branching
divisors $D_n$ resp. $\Delta_n$ are invariant under ${\tilde{G}}$ resp. $G$.

\vs

\noin {\bf 4.3. Lemma.} {\it A parametrized rational plane curve $C^* \s
\pr^{2*}$ is projectively equivalent to a parametrized quasi--monomial curve
iff $\pr^2_C \s \pr^n$ contains a one--dimensional $G$--orbit. This orbit is
unique iff $C^*$ is projectively equivalent to a parametrized quasi--monomial
curve, but not to a monomial one.}

\vs

\noin {\it Proof.} Let $\lambda \longmapsto (a_0 : \lambda a_1 :\dots :
\lambda^n a_n )$, where $\lambda \in \cz^*$, be a parametrization of the
$G$--orbit $O_p$ through the point $p =  (a_0 : \dots : a_n ) \in \pr^n$.
Since the non-zero coordinates are linearly independent as functions of
$\lambda$, the orbit $O_p \s \pr^n$ is contained in a projective plane iff all
but at most three of coordinates of $p$ vanish. If $p$ has exactly three
non--zero coordinates, then the only plane that contains ${\bar O}_p$ is a
coordinate one. If only two of the coordinates of $p$ are non--zero, then the
closure ${\bar O}_p$ is a coordinate axis. Since we consider a one--dimensional
orbit, the case of one non--zero coordinate is excluded. Now the lemma easily
follows from Lemma 4.2. \qed

\section{Proof of Theorem 3.1.$b''$}

\vs

In the sequel `bar' over a letter denotes a projective object, in contrast with
the affine ones.

\vs

\noin {\bf 5.1. Lemma.} {\it Let ${\bar H}_0$ be the hyperplane in $\pr^{n-1}$
given by the equation $\sum\limits_{i=1}^n x_i = 0$, and let ${\bar D}_{n-1} =
\bigcup\limits_{1\le i<j\le n} {\bar D}_{ij}$ be the union of the diagonal
hyperplanes, where ${\bar D}_{ij} \s  \pr^{n-1}$ is given by the equation
$x_i - x_j =0$. Then ${\bar H}_0 \se {\bar D}_{n-1}$ is C--hyperbolic.}

\vs

\noin {\it Proof.} Put $y_i = x_1 - x_{i+1}\,\,,\,\,i=1,\dots,n-1$. Then $z_i =
y_i / y_{n-1}\,\,,\,\,i=1,\dots,n-2 $, are coordinates in the affine chart
${\bar H}_0 \se {\bar D}_{1, n} \cong \cz^{n-2}$. In these coordinates ${\bar
D}_{1, i+1} \cap {\bar H}_0$ resp. ${\bar D}_{i+1, n} \cap {\bar H}_0$ is given
by the equation $z_i = 0$ resp. $z_i =1\,\,,\,\,i=1,\dots,n-2 $. Thus, ${\bar
H}_0 \se {\bar D}_{n-1} \hookrightarrow (\cz^{**} )^{n-2}$, where
$\,\cz^{**} := \pr^1 \se \{3\,\,points\}$. By Lemma 2.3 it follows that
${\bar H}_0 \se {\bar D}_{n-1}$ is C--hyperbolic. \qed

\vs

\noin {\it Remark.} Using the criterion
in Zaidenberg [35], (3.4) one can easily verify that
${\bar H}_0 \se {\bar D}_{n-1}$ is Kobayashi complete hyperbolic and
hyperbolically embedded into ${\bar H}_0 \cong \pr^{n-2}$.
Observe, by the way, that for $n = 4$, ${\bar D}_3 \cap {\bar H}_0$ is a
complete quadruple \footnote{i.e. a union of six lines through four points in
general position.} in ${\bar H}_0 \cong \pr^2$.

\vs

The proof of Theorem 3.1.$b''$ will be done in several steps.

\vs

\noin {\it Basic construction.} Let $q_0$ be a cusp of $C^*$ and
$q_0^* \s L_C \s \pr^2$ be the dual line. We will assume that
$q_0 = (0 : 0 : 1)$, so that $q_0^* = l_2 = \{x_2 = 0\}$. Due to Lemma 4.1,
in what follows we will
fix a special normalization $\nu\,:\,\pr^1 \to C^* \s \pr^{2*}$. Recall
(see sect. 4.a) that
$\nu (\infty ) = q_0$ and
$\pr^2_C = \rho_C (\pr^2 ) \s {\bar A}_1 \s \pr^n = S^n \pr^1$, where
${\bar A}_1 = \{(a_0 :\dots : a_n) \in \pr^n\,|\,a_1 = 0\} $. Set $\cz^n_a
= \{(a_1,\dots,a_n)\} = \{a \in \pr^n\,|\,a_0 \neq 0\},\,\, \cz^n_z =
s_n^{-1}(\cz^n_a) \s (\pr^1)^n,\,\, A_1 = {\bar A}_1 \cap \cz^n_a
\cong \cz^{n-1}$ and
$H_0 = \{z = (z_1,\dots, z_n ) \in \cz^n_z \,|\,
\sum\limits_{i=1}^n z_i =0 \} = s_n^{-1} (A_1)  \cong \cz^{n-1}$.
We have $\rho_C (X) \s \rho_C (\pr^2 \se l_2 ) \s \cz^n_a \s \pr^n$.
The restriction of the Vieta map yields the non--ramified covering
$s_n\,:\,{H}_0 \se D_n \to A_1 \se \Delta_n$.

Denote by $\pi$ the canonical projection $\cz^n_z \se \{0\} \to \pr^{n-1}$. Put
${\bar H}_0 := \pi ({H}_0 ) \cong \pr^{n-2} \s \pr^{n-1}$ and ${\bar D}_{ij} :=
\pi (D_{ij} )\,,\,{\bar D}_{n-1} := \pi (D_n ) = \bigcup\limits_{1\le i < j \le
n} {\bar D}_{ij}$. By Lemma 5.1
${\bar H}_0 \se {\bar D}_{n-1}$ is  C--hyperbolic.

Thus, we have the following commutative diagram:

\begin{picture}(250,80)
\unitlength0.2em
\thicklines
\put(34,5) {$X$}
\put(35,25){$Y$}
\put(47,25){$\hookrightarrow$}
\put(27,16){${\tilde s}_n$}
\put(36,22){{\vector(0,-1){11}}}
\put(47,30){${\tilde\rho}_C$}
\put(47,5){$\hookrightarrow$}
\put(47,11){$\rho_C$}
\put(62,5){$A_1 \se \Delta_n$}
\put(62,25){${H}_0 \se D_n$}
\put(70,22){\vector(0,-1){11}}
\put(75,16){$s_n$}
\put(170,15){(7)}
\put(82,27){{\vector(1,0){14}}}
\put(87,30){$\pi$}
\put(100,25){${\bar H}_0 \se {\bar D}_{n-1} \hookrightarrow \pr^{n-2}$}
\end{picture}

\noin where ${\tilde{s}}_n$ is the induced non--ramified covering.
  Note that the Vieta map $s_n$ is equivariant with respect to the
$\,\cz^*$--actions ${\tilde{G}}$ on $ {H}_0 \se D_n$ and $G$ on $A_1 \se
\Delta_n$, respectively, and all the fibres of the projection $\pi$ are
one--dimensional ${\tilde{G}}$--orbits (see (4.c)).

\vs

Since ${\bar H}_0 \se {\bar D}_{n-1}$ is C--hyperbolic (see Lemma 5.1),
by Lemma 2.3 we get that $Y$, and therefore also $X$, is almost
C--hyperbolic as soon as  all the fibres of the projection
$\pi \circ {\tilde \rho}_C \,:\,Y \to
{\bar H}_0 \se {\bar D}_{n-1}$ are finite. To prove Theorem 3.1.$b''$ it is
enough to check that this is the case under the condition ($b''$).

\vs

\noin {\it Claim 1. If the dual curve $C^*$ (as a parametrized curve)
is not projectively equivalent to a quasi--monomial one, then the mapping $\pi
\circ {\tilde \rho}_C \,:\,Y \to {\bar H}_0 \se {\bar D}_{n-1}$ has finite
fibres.}

\vs

\noin Indeed, since the fibres of $\pi$ are $\tilde G$--orbits, it is enough to
show that any $\tilde G$--orbit in $H_0 \s \cz^n_{(z)}$ has a finite
intersection with ${\tilde \rho}_C (Y)$. Or, what is equivalent, that any
$G$--orbit in $A_1 \s \cz^n_{(a)}$ has a finite intersection with
$\rho_C (X) \s \pr^2_C$. We have shown in Lemma 4.3 above that if the latter
fails, i.e. if $\pr^2_C$ contains a one-dimensional $G$--orbit, then $C^*$
(parametrized as above) is projectively equivalent to a (parametrized)
quasi--monomial curve, which is assumed not to be the case. This yields claim
1. \qed

\vs

Thus, we may suppose that $C^*$ (as a parametrized curve) is projectively
equivalent
to a quasi--monomial one. By Lemmas 4.2, 4.3 this means that $\pr^2_C$ contains
a coordinate line, which is the closure of a one-dimensional $G$--orbit $O_p$.
Any monomial curve of degree $n$ with a cusp is
projectively equivalent
(as a parametrized curve, in appropriate parametrization) to one of the curves
$(1 : t^k : t^n)$, where $1 \le k \le n-2$. But this is excluded by the
conditions in ($b''$). Hence,
$C^*$ (as a parametrized curve, with the special normalization chosen above)
can not be  projectively equivalent to a monomial curve, i.e. the plane
$\pr^2_C$ is not a coordinate one (see Lemma 4.2).

\noin Let $l_{n-k,\,n-m} \s \pr^2_C\,,\,\,0\le k < m \le n$, be the only
coordinate axis contained in $\pr^2_C$.
We will distinguish between two cases:

\ms

\noin (i) $\,\,\,l := \rho_C^{-1} (l_{n-k,\,n-m}) \s L_C$ {\hskip 0.5in}
and
{\hskip 0.5in} (ii) $\,\,\,l = \rho_C^{-1} (l_{n-k,\,n-m})
\not\subset L_C\,\,.$

\ms

\noin Note that (i) resp. (ii) holds iff the dual point $q = l^* \in \pr^{2*}$
is resp. is not a cusp of $C^*$. \\
If (i) holds then, as in  Claim 1 above, $\pi \circ
{\tilde \rho}_C \,:\,Y \to {\bar H}_0 \se {\bar D}_{n-1}$ has finite fibres,
and hence $X$ is almost C--hyperbolic. \\
Thus, the following claim finishes the proof of ($b''$).

\vs

\noin {\it Claim 2.
Assume that the dual curve $C^*$ (as a parametrized curve)
is projectively equivalent to a quasi--monomial, but not to a monomial one.
Then (ii) holds iff we have one of the two exceptional cases
in ($b''$).}

\vs

\noin {\it Proof.}
Let $\nu = (h_0 : h_1 : t^n + h_2)$, where $h_i \in \cz[t]$ and
deg$\,h_i \le n-2\,,\,\,i=0,1,2$, be the special normalization of $C^*$
fixed in the basic construction above. The inclusion $l_{n-k,\,n-m} =
\,$span$\,(e_{n-k},\,e_{n-m}) \s \pr^2_C$ means that
$t^k,\,t^m \in \,$span$\,(h_0,\,h_1,\,t^n + h_2)$. Consider two cases

\ss

\noin (1) $k < m \le n-1\,,$ {\hskip 0.5in} and
{\hskip 0.5in} (2) $k < m = n\,\,.$

\ss

\noin In the first case we have $t^k,\,t^m \in \,$span$\,(h_0,\,h_1)$, so
that $k,\,m \le n-2$, and without lost of generality we may suppose that
$\nu = (t^k : t^m : t^n + h_2(t))$. In the second case we have $t^k,\,h_2 \in
\,$span$\,(h_0,\,h_1)$ (so, in particular, $k \le n-2$), and we may assume
that $\nu = (t^k : h_1(t) : t^n)$. With these conventions we have
$B_C (e_0) = e_{n-k}$; $B_C (e_1) = e_{n-m}$ in the case (1), so that $l =
\rho_C^{-1} (l_{n-k,\,n-m}) = l_2$, and
$B_C (e_2) = e_{0}$ in the case
(2), so that $l = \rho_C^{-1} (l_{n-k,\,n-m}) = l_1$. \\
Thus, in the first case $l^* = l_2^* = q_0 = (0 : 0 : 1)$ is the cusp of
$C^*$, and hence (i) holds. In the second case the dual point
$l^* = l_1^* = q_1 = (0 : 1 : 0)$ is a cusp of $C^*$ iff $k \ge 2$ (
indeed, $q_1 \in C^*$ is a smooth point if $k=1$, and
$q_1 \notin C^*$ if $k=0$).
Thus, (ii) holds iff $\nu = (t^k : h_1(t) : t^n)$ with $k \le 1$ and
deg$\,h_1 \le n-2$, which are exactly the two exceptional cases of $(b'')$.
This completes the proof of Claim 2 and hence the proof of $(b'')$. \qed

\vs

\noin {\it Remarks.}
1. If $C^*$ is one of the exceptional curves mentioned in (3.1.$b''$) and is
not projectively equivalent to a monomial curve, then
$X = \pr^2 \se (C \cup L_C)$ is C--hyperbolic modulo $l_1 = \rho_C^{-1}
(l_{n-k,\,0})$ (in a natural sense). But this is not true, in general, for a
plane curve whose dual is a quasi-monomial curve without cusps. An example
is a three--cuspidal plane quartic $C \s \pr^2$. Indeed, then $C^*$ is a nodal
cubic,
which is projectively equivalent to a quasi--monomial curve $t \longmapsto
 (t : t^3 : t^2 - 1 )$.
The Kobayashi pseudo--distance of $\pr^2 \se C$ is degenerate on at
least seven lines \footnote{these lines are: the three cuspidal tangents, the
three lines through a pair of cusps and the only bitangent line.}, and thus
$\pr^2 \se C$ is not C--hyperbolic modulo
a line. Furthermore, $\pi_1(\pr^2 \se C)$ is a finite non-abelian group
of order $12$ (see Zariski [38], p.143). Thus, any covering over $\pr^2 \se C$ is
a Liouville one.

\vs

\noin 2. Let $C^*$ be a monomial curve, and both $C^*$ and $C$ belong to
the linear pencil
$C_{\mu} = \{\alpha x_0^n + \beta x_1^k x_2^{n-k} =0\}$, where
$\mu = (\alpha : \beta ) \in \pr^1$. Set $X = \pr^2 \se (C \cup L_C)$. Then
the Kobayashi pseudodistance $k_X$ is degenerate along any of the members of
this linear pencil. At the same time, the distance between points on two
distinct members is
always positive. In particular, any entire curve $f\,:\,\cz \to X$ is
contained in
one of the curves $C_{\mu}$.

\vs

The above proof gives us an additional information that will be used in
the concrete examples of Sect. 6.$a$ to distinguish the exceptional cases. It
can be summarized as follows.

\vs

\noin {\bf 5.2. Proposition. } {\it Let
$\nu = (h_0 : h_1 : t^n + h_2)\,:\, \pr^1 \r C^*$ be a special
normalization. Then $X =\pr^2 \se (C \cup L_C)$ is
C-hyperbolic in each of the following cases:\\
a) $C^*$ has at least three cusps.\\
b) $C^*$ has two cusps and span$\,(1,\,h_2) \not\subset$ span$\,(h_0,\,h_1)$.\\
c) $C^*$ has one cusp, span$\,(1,\,h_2) \not\subset$ span$\,(h_0,\,h_1)$, and
span$\,(t,\,h_2) \not\subset$ span$\,(h_0,\,h_1)$.}

\vs

The proof is easy and can be omited.

\section{Examples}
Hereafter ${\cal H}(d)$ denotes the set of all plane curves of degree
$d$ with complete hyperbolic and hyperbolically embedded complements.
In this section we give explicit examples of plane curves with C--hyperbolic
complements. Furthermore, we construct, for every even $d \geq 6$, families
of irreducible curves in
${\cal H}(d)$, especially of elliptic or rational such curves.
They are described by the degree, the genus and the singularities of their
members or of the dual curves. Most of them arise from Theorem 1.1 in the
special case where $C^*$ is a nodal curve, and only the case of maximal
cuspidal rational sextics has to be treated in a different way (see
Proposition 6.11 below).

\vs

\noin {\bf a) Reducible curves}

\vs

\noin {\bf 6.1. Two examples of quintics}

\ss

\noin 1. Perhaps, the simplest example is
the arrangement $C_5$ of five lines with two triple points. It is projectively
unique and can be given by the equation $x_0 x_1 x_2 (x_0 -x_1 )(x_0 -x_2)
=0\,\,.$ The complement $X = \pr^2 \se C_5$ is biholomorphic to $(\cz^{**})^2$,
and thus its universal covering is the bidisc $U^2$. Hence, $X$ is C-hyperbolic
and also complete hyperbolic. However, by  Lemma 2.2 $X$ is not hyperbolically
embedded into $\pr^2$.

\vs

\noin 2. Another example is a smooth conic $C$ together with its three distinct
tangents $L = l_1 \cup l_2 \cup l_3$. This configuration is also projectively
unique. We may identify $C$ with the discriminant $\Delta_2 \s \pr^2$, so that
the Vieta covering $s_2\,:\,(\pr_1)^2 \to \pr^2$ is branched along $C$, which
is the
image of the diagonal $D_2 \s (\pr^1)^2$. The lines $l_1,\,l_2,\,l_3$ are the
images of six generators of the quadric $\pr^1 \times \pr^1$, three horizontal
ones and three vertical ones.
Thus, $X = \pr^2 \se (C \cup L)$ is covered by $(\cz^{**})^2 \se D_2$,
being, henceforth, C-hyperbolic.

\vs

\noin {\bf 6.2. Four examples of sextics}

\ss

\noin 1. Modifying example 6.1.1 consider
an arrangement $C_6$ of six lines with three triple points. This is a
complete quadruple (cf. the remark after Lemma 5.1). It is
projectively unique and can be given by the equation
$x_0 x_1 x_2 (x_0 -x_1 )(x_0 -x_2)(x_1 - x_2) =0\,\,.$
It is known (see e.g. Kaliman [19]) that the universal covering of the complement
$X = \pr^2 \se C_6$ is biholomorphic to the Teichm\"uller space
$T_{0,\, 5}$ of the Riemann sphere with five punctures. Thus, via the Bers
embedding $T_{0, \,5} \hookrightarrow \cz^2$ it is biholomorphic to a bounded
Bergman domain of holomorphy in  $\cz^2$, which is contractible and Kobayashi
complete hyperbolic. The automorphism group of $T_{0,\,5}$ is discrete and
isomorphic to the mapping class group, or modular group, ${\rm Mod}(0,\,5)$
(see Royden [30]). Clearly, the fundamental group $\pi_1(X)$ is a subgroup
of finite index in ${\rm Mod}(0,\,5)$.

\vs

\noin 2. The next three examples serve as illustrations to ($b''$) of Theorem
3.1. Consider a nodal cubic $C \s \pr^2$ together with its three  inflexional
tangents $L_C = l_1 \cup l_2 \cup l_3$.
They correspond to the cusps of the dual curve $C^*$, which is a 3-cuspidal
quartic. Both $C$ and $C^*$ are projectively unique. By Proposition 5.2.$a$, we
have that
$X = \pr^2 \se (C \cup L_C)$ is almost C--hyperbolic.

\vs

\noin 3.
Let $C \s \pr^2$ be the rational quintic
$t \longmapsto (2t^5 - t^2 : -(4t^3 + 1) : 2t)$ with a cusp at the only
singular point $(1 : 0 : 0)$. The dual curve $C^*$ is the quasi--monomial
quartic $t \longmapsto (1 : t^2 : t^4 + t)$ given by the equation $(y_0 y_2 -
y_1^2 )^2 = y_0^3 y_1$. It has the only singular
point $q_0 = (0 : 0 : 1)$, which is a ramphoid cusp, i.e. it has the
multiplicity sequence\footnote{Recall that  {\it the multiplicity sequence}
of a plane analytic germ $A$ at $p_0 \in A$ is the sequence of multiplicities
of $A$ at $p_0$ and in its infinitesimaly near points.}
$(2,\,2,\,2,\,1,\,\dots)$
and $\delta = \mu /2 = 3$,
where $\mu$ is the Milnor number. Any rational quartic with a ramphoid cusp
is projectively equivalent to $C^*$ (see Namba [25], 2.2.5(a)). The artifacts $L_C$
consist of the only cuspidal tangent line $l_2 = \{x_2 = 0 \}$  of $C$. By
Proposition 5.2.$c$, the complement $X = \pr^2 \se (C \cup l_2)$ is almost
C--hyperbolic.

Note that a smooth affine curve $\Gamma = C \se l_2 \s  \pr^2 \se l_2
\cong \cz^2$ is isomorphic to $\cz^* := \cz \se \{0\}$, and so its
complement $X := \cz^2 \se \Gamma$ is almost C--hyperbolic.

\vs

\noin 4. Let $C' \s \pr^2$ be the rational quartic
$t \longmapsto (t^3 (2t + 1)  : -t(4t + 3) : -2)$.
It has two singular points, a double cusp at the point
$(0 : 0 : 1)$ (i.e. a cusp with the multiplicity sequence
$(2, \,2,\,1,\,\dots)$ and $\delta = 2$) and another one, which is an
ordinary cusp. The dual curve   $C'^* \s \pr^{2*}$ is the quasi--monomial
quartic $t \longmapsto (1 : t^2 : t^4 + t^3)$ given by the equation
$(y_0 y_2 - y_1^2 )^2 = y_0 y_1^3$. It has the same type of singularities
as $C'$, namely a double cusp at the point $q_0 = (0 : 0 : 1)$ and an
ordinary cusp at the point $(1 : 0 : 0)$. Therefore, $L_{C'} = l_0 \cup l_2$,
where $l_0 = \{x_0 = 0 \}$ and $l_2 = \{x_2 = 0 \}$. Put $s = 1/t$ and permute
the coordinates to obtain the special
normalization $s \longmapsto (1+s:s^2:s^4)$ of $C'^*$. Now by Proposition
5.2.$b$,
the complement $X:= \pr^2 \se (C' \cup L_{C'})$ is almost C--hyperbolic.

\vs

\noin {\bf 6.3. Example of a septic}

\ss

\noin Let things be as in example 6.1.2. Performing the Cremona transformation
$\sigma$ of $\pr^2$ with center at the points of intersections of the lines
$l_1,\, l_2,\, l_3$, we obtain a 3-cuspidal quartic $C' :=  \sigma (C)$
together with three new lines $L' = m_1 \cup m_2 \cup m_3$ passing through
pairs of cusps of $C'$. Put $X' = \pr^2 \se (C' \cup L')$. Since
$X = \pr^2 \se (C \cup L)$ is C--hyperbolic and $\sigma\,|\,X\,:\,X \to X'$
is an isomorphism, $X'$ is also C--hyperbolic.

\vs

\noin {\bf 6.4. Two examples of octics}

\ss

\noin Next we pass to examples to part $a$) of Theorem 3.1. Let $C^* \s
\pr^{2*}$ be an irreducible Pl\"ucker curve of
genus $g\ge 1$ with $\kappa$ cusps. Then the dual curve $C \s \pr^2$ has $f =
\kappa$ ordinary flexes, and $L_C$ is the union of inflexional tangents of $C$.
By the class formula (1), we have
$d = {\rm deg}\,C = 2(n+g-1)- \kappa$. Since all $\kappa$ inflexional
tangents of $C$ are distinct, it follows that ${\rm deg}\,(C\cup L_C) =
2(n+g-1) \ge 2n \ge 6$. Assume that $\kappa>0$ to exclude the case when $C^*$
is an immersed curve (cf. (6.5), (6.6) below). Since $g\ge 1$, the case
when $C$ is a singular cubic has also been excluded. Thus, we have $n \ge 4$,
and hence ${\rm deg}\,(C\cup L_C) \ge 8$.

\ss

\noin 1. The simplest example is a quartic $C^*$ with an ordinary cusp
and a node
as the only singularities (see Namba [25], p.133).
The dual curve $C$ is an elliptic septic with the only inflexional tangent
line $l=L_C$.

\ss

\noin 2. Another example is a quartic $C^*$ with two ordinary cusps as
the only
singular points (see Namba [25], p.133). Here $C$ is an elliptic
sextic and $L_C$ is the union of two inflexional tangents of $C$.

\ss

\noin In both examples the assumptions of Theorem 3.1.$a$ are fulfilled,
and so $X = \pr^2 \se (C \cup L_C)$ is C--hyperbolic.

\vs

\noin {\it Remark.} It can be checked that in examples
6.1.2, 6.2.1, 6.2.2, 6.3, 6.4.1 and 6.4.2 the conditions of Lemma 2.2 are
fulfilled, and hence the
corresponding complements are Kobayashi complete hyperbolic and
hyperbolically embedded into $\pr^2$, whereas in 6.1.1, 6.2.3 and 6.2.4
hyperbolic embeddedness fails.

\vs

\noin {\bf b) Irreducible curves}

\vs

\noin {\bf 6.5. Examples of irreducible curves of genus $g \ge 2$}

\ss

\noin
Theorem 1.1.$a$ can be applied, for instance, to an irreducible curve
$C \s \pr^2$ of genus $g \ge 2$ whose
dual $C^*$ is a nodal curve of degree $n \ge 4$ with $\delta$ nodes.
Such a curve $C$ does exist for any given $\delta$ with $0\le\delta
\le {n-1 \choose 2} -2$ (see Severi [31], Sect. 11, p.347; Oka [26], (6.7)).
By the class formula (1) and the genus formula $C$ has
degree $d=n(n-1) -2\delta$, which can be any even integer from the interval
$[2(n+1), n(n-1)]$. The minimal value of $d$ is $d = 10$, which corresponds
to a nodal quartic $C^*$ with one node (see Namba [25], p.130).

\vs

\noin {\it Remark.}
It was shown by Green [16], Carlson and Green [4] and Grauert and
Peternell [15] that an irreducible plane
curve $C$ of genus $g \ge 2$ belongs to ${\cal H}(d)$ if the following
conditions hold:

\noin (i) each tangent line to $C^*$ intersects with $C^*$ in at least
two points, and

\noin (ii) $2n < d$, where as before $d = {\rm deg}\,C$ and
$n = {\rm deg}\,C^*$.

These conditions are less restrictive than those above,
since here $C^*$ may possess cusps. For such a $C^*$  by the genus
formula
$2g \leq (n-1)(n-2)$, hence $n \geq 4$ for $g \geq 2$, and by (ii)
we have $d \geq 9$. Due to the class formula (1), this lower bound
is really achieved
for the family of duals of the irreducible quartics $C^*$ with an
ordinary cusp as the only singular point (see Namba [25], p.130).
However, we do not know whether in this example the complement of
$C$ is also C-hyperbolic.

\vs

\noin {\bf 6.6. Examples of elliptic curves}

\ss

\noin If the dual $C^*$ of $C$ is an immersed elliptic curve, then  by
the class formula (1) $d= {\rm deg}\,C =2n \ge 6$, where
$n = {\rm deg}\,C^* \ge 3$. Let $C$ be a sextic in $\pr^2$ with nine cusps.
Then $C$ is an elliptic Pl\"ucker curve whose dual $C^*$ is a smooth cubic;
vice versa, the dual of a smooth cubic is a sextic with nine ordinary cusps.
{}From Theorem 1.1.$a$ we get the following

\vs

\noin {\bf 6.7. Proposition.} {\it Every irreducible plane sextic with nine
cusps has C-hyperbolic complement and belongs to ${\cal H}(6)$.}

\vs

\noin Note that up to projective equivalence this family is one
dimensional. We refer to Gelfand, Kapranov and Zelevinsky [14], I.2.E for explicit
Schl\"afli's equations of these elliptic sextics. For instance, the dual of
the Fermat cubic $-x_0^3 + x_1^3 + x_2^3 = 0$ is the sextic
$$x_0^6 + x_1^6 + x_2^6 - 2x_0^3 x_1^3 - 2x_1^3 x_2^3 - 2x_0^3 x_2^3 =
0 \,\,.$$

Another example in degree $8$ is the family of elliptic curves
dual to the nodal quartics with two nodes (see e.g. Namba [25], p.133 for the
existence). Together with (6.4) and (6.5) this yields the following

\vs

\noin {\bf 6.8. Proposition.} {\it For any even $d \ge 6$ there exists a
family of irreducible plane curves of degree $d$ and of genus $g \ge 1$
with C-hyperbolic complements which  belong to ${\cal H}(d)$. It is the
family of dual curves to the nodal Pl\"ucker curves of degree $n \ge 3$ with
$\delta$ nodes, with appropriate $n$ and $\delta$.}

\vs

\noin {\bf 6.9. Examples of rational curves}

\ss

\noin They illustrate ($b$) of Theorem 1.1.
A generic rational curve $C^*$ of degree $n \geq 3$
is a nodal Pl\"ucker curve (see e.g. Aure [2]). Its dual curve $C$ has an even
degree $d=2(n-1)$ and $\kappa = 3(d-2)/2$ cusps. Vice versa, any rational
Pl\"ucker curve $C$ of even degree $d=2(n-1)$ with $\kappa = 3(n - 2)$ cusps
is dual to a nodal curve $C^*$ of degree $n$. Here $\kappa$ is the
maximal number of
cusps which a rational Pl\"ucker curve of degree $d$ can possess,
and so these curves
are called {\it rational maximal cuspidal curves} (see Zariski [38], p. 267). Applying
Theorem 3.1.$b'$ we obtain the following

\vs

\noin {\bf 6.10. Proposition.} {\it For any even degree $d \geq 8$
a rational maximal cuspidal plane curve of degree $d$ belongs to
${\cal H}(d)$ and its complement is almost C-hyperbolic.}

\vs

What happens with rational maximal cuspidal curves of lower degrees?
For $d=4$ we have a three cuspidal quartic (which is projectively unique
(see Namba [25], p.146)). As we saw in the Remark 1 after the proof of Theorem 3.1.$b''$
 its complements is not even Kobayashi hyperbolic.
It remains the case $d=6$. In this case we have the following

\vs

\noin {\bf 6.11. Proposition.} {\it A generic rational maximal cuspidal
plane sextic belongs to ${\cal H}(6)$.}

\vs

\noin {\it Proof.} We keep the notation of sect. 2.b. From the proof of
Theorem 3.1.$b'$ we know that such a sextic $C$ is a generic plane
section of the
discriminant hypersurface $\Delta_4 \s \pr^4$ (by the plane
$\pr_C^2 = \rho_C (\pr^2)$). Clearly, being generic, $\pr_C^2$
does not meet the only
one--dimensional PGL$(2,\cz)$--orbit $O_4$. From the definition of the Zariski
embedding
it easily follows that it intersects the orbit closure
$S_1 = {\bar O}_{3,\,1}$ resp. $S_2 := {\bar O}_{2,\,2}$ in the set
$K=\{$the cusps of $C\}$ resp. $N=\{$the nodes of $C\}$. Therefore,
it intersects the only 3-dimensional orbit $O_{2,\,1,\,1}$ contained
in $\Delta_4$ in the curve $C \se (K \cup N)$.
 By the Pl\"ucker formulas, ${\rm card}\,K = 6$ and ${\rm card}\,N = 4$
 (this agrees with the fact that ${\rm deg}\,S_1 = 6$ and ${\rm deg}\,S_2 = 4$
(see Aluffi and Faber [1], Proposition 1.1)). Let $C_q = \pr_C^2 \cap {\bar O}_q$,
 where $O_q = O_{1,\,1,\,1,\,1}$, i.e. $q \in \pr^4 \se \Delta_4$.
 Since ${\bar O}_q = O_q \cup S_1$ (see sect. 2.b) it is easily seen that
 the curve $C_q \s \pr^2_C$ meets $C$ exactly in the cusps of $C$.

Now we use diagram (6). Let $f: \cz \r \pr_C^2 \se C = \pr_C^2
\se \Delta_4 $ be an entire curve. Since $s_4\,:\,(\pr^1)^4 \se D_4 \r
\pr^4 \se \Delta_4$ is an unramified covering, $f$ can be lifted to $(\pr^1)^4
\se D_4$. The curve $\,\cz^{**}$ being hyperbolic, this lifted entire curve
has to be contained in a fiber of $\tilde{\pi}_4$, which is an orbit of
the PGL$(2,\cz)$-action on $(\pr^1)^n$. The Vieta map $s_4$ being equivariant,
the entire curve $f\,:\, \cz \r \pr_C^2 \se C$ is contained
in a PGL$(2,\cz)$--orbit, too.

Thus, to see that $\pr^2 \se C$ is Brody hyperbolic it is enough to show
that the quasiprojective curves $C_q \se C = O_q \cap \pr^2_C$ are hyperbolic
for all $q \in \pr^4 \se \Delta_4$. Once this is done, Proposition 6.11
follows from Lemma 2.2.

It is well known (see Hilbert [18], p.58 or Popov and Vinberg [29]) that the 3-dimensional
PGL$(2,\cz)$--orbit
closures in $\pr^4$ form a linear pencil. This pencil of sextic threefolds
is generated by its members $3P$ and $2H$, where the irreducible quadric resp.
cubic $P$ and $H$ are defined by the basic invariants
$\tau_2 = a_0 a_4 - 4a_1 a_3 + 3a_2^2$ resp.
$\tau_3 = a_0 a_2 a_4 - a_0 a_3^2 - a_1^2 a_4 + 2a_1 a_2 a_3 - a_2^3$
(here we use the coordinates where $q(u,\,v) = a_0 u^4 + 4a_1 u^3 v +
6a_2 u^2 v^2 + 4a_3 uv^3 + a_4 v^4$). The base point set of this linear
pencil is the surface $S_1 = H \cap P$, as it follows from the description
of the orbit closures in Aluffi and Faber [1] (see sect. 2.b).

The restriction of the above pencil to the plane $\pr^2_C$ is the linear
pencil of plane sextics $\alpha := (C_q)$ generated by $3p$ and $2h$,
where $p := P \cap \pr^2_C$ and $h := H \cap \pr^2_C$ are respectively
irreducible conic and cubic. Its base point set $p \cap h$ is the set $K$
of cusps of $C$ (note that $C$ itself is a member of $\alpha$). The
intersection of $p$ and $h$ at the points of $K$ is transversal,
because ${\rm card}\,K = p\cdot h = 6$. Since the ideal generated by two
distinct members $C' = C_{q'}$ and $C'' = C_{q''}$ is the same as the one
generated by $3p$ and $2h$, we have for the local intersection multiplicities
at any point $x \in K$
$$ i(C',\,C'';\,x) = i(3p,\,2h;\,x) = 6\, i(p,\,h;\,x) = 6\,.$$

Assume now that a member $C_q$ of the pencil $\alpha$ has an irreducible
component $T$ which intersects $C$ in at most two points $x', \,x'' \in K$.
Since $i(T,\, C;\,x) \leq 6$ for $x = x' ,\,x''$, we would have
deg$\,C \cdot$deg$\,T \leq 12$, and hence deg$\,T \leq 2$.

If $T$ would be
a projective line, then by Bezout's Theorem $T\cdot C = {\rm deg}\,C = 6$,
hence $i(T,\, C;\,x) = 3$ for $x = x' ,\,x''$, and so $T$ should be a
common cuspidal tangent of $C$ at these two cusps $x' ,\,x''$, which is
impossible for a Pl\"ucker curve $C$.

If, further, $T$ would be a smooth
conic, then by Bezout's Theorem we would have $i(T,\, C;\,x) = 6$ for
$x = x' ,\,x''$, again in contradiction with the fact that $C$ is a
Pl\"ucker curve. Indeed, an ordinary cusp $(C,\,P)$ can be uniformized
by $t \longmapsto (t^2, t^3 + O(t^4))$, see e.g. Namba [25], 1.5.8,
and therefore the local intersection multiplicity of an ordinary cusp with
a smooth curve germ $(C',\,P)$ can be at most 3. To see this, observe that
plugging this parametrization into the power series expansion at $P$ of the
defining equation of $C'$, its linear term will contain a non--zero
monomial in $t$ of order at most $3$, which cannot be cancelled by further
higher order terms \footnote{see also Fulton [13], (1.4) for a more general
fact.}.

Thus, there is no irreducible component $T$ as above, and hence all
the non--compact curves $C_q \se C$ in $\pr_C^2 \se C$ are hyperbolic.
\qed

\vs

\noin {\it Remark.} The dimension of the family of all plane rational
nodal curves of degree $n \ge 3$ modulo projective equivalence is
$3(n-3) = \frac{3}{2}(d-4)$, where $d=2(n-1)$ is the degree of the dual curves.
In particular, the family of curves in Proposition 6.11 is three--dimensional.

\section{Miscellaneous}

\noin {\bf a) Plane curves with a big fundamental group of the  complement}

\vs

\noin
Due to Lin's Theorem mentioned in the Introduction (see Lin [22], Thm. 13),
we obtain the following

\vs

\noin {\bf 7.1. Proposition.} {\it If $C \s \pr^2$ is one of the curves
mentioned in Theorem 3.1, then the group $\pi_1( \pr^2 \se (C \cup L_C))$
is not almost nilpotent. In particular, it is so in all
the examples of sect. 6.}

\vs

\noin Note that in certain cases more strong fact holds. Let us say that
a group $G$ is {\it big} if it contains a non--abelian free subgroup.
By a theorem of von Neumann, a big group is non--amenable. The converse
is not true in general; the corresponding
examples are due to  A. Ol'shanskiy, S. I. Adian and M. Gromov  (see e.g.
Ol'shanskiy and Shmel'kin [27]). But the groups $G$ in all these examples
are not finitely presented. For finitely presented
groups the equivalence of bigness and non--amenability is unknown
\footnote{we are thankful to V. Sergiescu and V. Guba for this information.}.
Being non--amenable, a
big group can not be almost nilpotent or even almost
solvable. As follows from the Nielsen--Schreier Theorem, a subgroup of
finite index of a big group is big, as well
as a normal subgroup with a solvable quotient. Clearly, a group with
a big quotient is big.

The following conjecture seems to be plausible.

\vs

\noin {\bf 7.2. Conjecture.} {\it If an algebraic variety $X$ is
C--hyperbolic, then $\pi_1(X)$ is a big group.}

\vs

Note that by another Lin's theorem (Lin [22], Thm. B(b)), $\pi_1(X)$ as above
can not be an amenable group with a non--trivial center, at least if the
universal covering space $\tilde X$ is Carath\'eodory hyperbolic. Observe
also that the conjecture is obviously true for dim$\,X=1$.

\ss

As far as the complements of plane curves is concerned, we have the following
fact\footnote{its proof, which was done jointly with S. Orevkov, will be
published elsewhere.}.

\vs

\noin {\bf 7.3. Theorem.} {\it Let $C \s \pr^2$ be an irreducible curve whose
dual $C^*$ is an immersed curve of degree $n$ and of geometric genus $g$,
where $n \ge 2g+1$ and $n \ge 4$ if $g = 0$. Then the group
$\pi_1(\pr^2 \se C)$ is big.}

\vs

\noin {\it Remarks.} 1. We do not know whether the theorem
is true without the assumption $n \ge 2g+1$.

\noin 2. A presentation of the group $\pi_1(\pr^2 \se C)$
for a generic maximal cuspidal curve $C \s \pr^2$ of genus  0 or 1
was found by Zariski
[38], p. 307; c.f. also Kaneko [20] for the case $n \geq 2g+1$, where $n = {\rm
deg}\,C^*$.
However, even if such a presentation is given it might be not so easy to
deduce Theorem 7.3.
\vs

{\bf b) Minimal degree of an irreducible curve with C--hyperbolic complement}

\vs

Here we show that the examples in sect. 6.b are,
indeed, at the borderline, as far as the C-hyperbolicity is concerned.
Observe that for curves of degree $\leq 4$ the
complement is not even hyperbolic, since there always exist projective
lines which intersect $C$ at most in two points, see e.g. Green [17].
The same remains true for irreducible quintics which are not Pl\"ucker.

\vs

\noin {\bf 7.4. Lemma.} {\it Let $C \s \pr^2$ be an irreducible quintic
which is not a Pl\"ucker curve. Then $\pr^2 \se C$ is not Brody hyperbolic.
Moreover, there exists a line $l_0 \s \pr^2$ which intersects with $C$ in at
most two points.}

\vs

\noin {\it Proof.} Assume that $C$ has a non--classical singular point $p_0$
(see sect. 2.a). Let $l_0$ be the tangent line to a local analytic branch of
$C$ at $p_0$. If ${\rm mult}_{p_0}\,C \ge 3$, then $i(C,\,l_0 ;\,p_0 ) \ge 4$,
and so $l_0$ intersects with $C$ in at most one more point.  If ${\rm
mult}_{p_0}\,C = 2$, then either  $p_0 \in C$ is a tacnode, i.e. $C$ has two
smooth branches at $p_0$ with the same tangent $l_0$, or
$C$ is locally irreducible in $p_0$ and has the multiplicity sequence
$(2,\,2,\,\dots)$ at $p_0$.
In both cases we still have  $i(C,\,l_0 ;\,p_0 ) \ge 4$, and the same
conclusion as before holds. It holds also in the case when $l_0$ is the
inflexional tangent to $C$ at a point where $C$ has a flex of order at
least $2$ (see Namba [25], (1.5)).

Therefore, we may suppose that $C$ has only
classical singularities and ordinary flexes.
Let $q_0$  be a singular point of $C^*$ which is not classical.
It can not be locally irreducible, since $C$ has only ordinary flexes.
If one of the local branches of $C^*$ at $q_0$ is singular,
then the dual line $l_0$ of $q_0$ is an inflexional tangent at some flex
of $C$, tangent also at some other point. By Bezout's Theorem
$l_0$ is a bitangent line with intersection indices $2$ and $3$.

It remains to consider the case when $C^*$ has only smooth local
branches at $q_0$. If two of them, say, $A^*_0$ and $A^*_1$, are tangent
to each other, then by duality the corresponding local branches $A_0$
and $A_1$ of $C$ should have common center and also be
tangent to each other. This is impossible since $C$ is supposed
to have only classical singularities. Thus,
$q_0 \in C^*$ should be an ordinary singular point with at least three
distinct
branches. But then the dual line $l_0$ of $q_0$ is tangent
in at least three different points of $C$, which contradicts to
Bezout's Theorem. \qed

\vs

\noin {\bf 7.5. Proposition.} {\it The minimal possible degree of
an irreducible plane curve with
almost C-hyperbolic complement is six, and
it has singularities  worse than ordinary double points.}

\vs

\noin {\it Proof.} By Lemma 7.4 and the preceding remarks, to prove the
inequality $d \ge 6$ it is enough to exclude the Pl\"ucker quintics.
Due to Degtyarev's list [5,6],
the fundamental group of the complement of an irreducible Pl\"ucker quintic
is abelian, and so it is isomorphic to $\gz / 5 \gz$. Thus, the only
non--trivial
covering $Y$ over $\pr^2 \se C$ is a finite cyclic one. Being
quasiprojective, such an $Y$ is a Liouville variety, and hence
$\pr^2 \se C$ is not almost C--hyperbolic.

\noin The second statement follows from the theorems of Deligne--Fulton
and Lin, as it was explained in the introduction. \qed

\vs

\noin {\bf c) Genericity of the inflexional tangent lines}

\vs

In Theorem 3.1 we gave sufficient conditions of
(almost) C-hyperbolicity of the
complement of a plane curve together with its artifacts. Since the
complement of the curve itself is only rarely C--hyperbolic (in particular,
this never happens for a nodal curve, see the discussion in the introduction),
in order to guarantee C-hyperbolicity
we need to add the artifacts, or  at least some of them.
But then the question arises
whether the complement to artifacts themselves is C-hyperbolic.
This is the case, for instance, when the artifacts contain the
configuration of five lines as in example 6.1. Our aim here is to show
that this is not the case for a generic plane curve. Observe that
being generic such a curve is smooth, and hence its artifacts
are just the inflection tangents.

\vs

\noindent {\bf 7.6. Proposition.}
{\it If $C$ is a generic plane curve, then the artifacts $L_C$
are in general position.
In particular, $\pr^2 \se L_C$ is not C-hyperbolic.}

\vs

\noin {\it Proof.} Since a generic smooth
curve is a Pl\"{u}cker curve, we know that all its inflection tangents
are distinct.
Consider the quasiprojective variety ${\cal L}$ of all the configurations
$l=(l_1,l_2,l_3,P_1,P_2,P_3)$,
where $l_1$, $l_2$, $l_3$ are three distinct lines in $\pr^2$ all passing
through a common point, and $P_i \in l_i,\,i=1,2,3$, are pairwise distinct
points. Let
${\cal C}(d)$ be the quasiprojective variety of the smooth plane curves
of degree $d$. For a given $l=(l_1,l_2,l_3,P_1,P_2,P_3) \in {\cal L}$
denote by ${\cal C}(d)_l$ the subvariety of curves in ${\cal C}(d)$
which have flexes at $P_i$ with inflectional tangents $l_i,\,i=1,2,3$.

\vs

\noin {\it Claim. For any $d \geq 4$ and for all $l \in {\cal L}$
we have} codim$\,_{{\cal C}(d)} {\cal C}(d)_l = 9$.
\vs

\noin {\it Proof.}
It can be easily shown that there are exactly three Aut$(\pr^2)$--orbits in
${\cal L}$, say ${\cal L}_a,\,{\cal L}_b,\,{\cal L}_c$,
where ${\cal L}_a$ is the only open orbit which consists of the
configurations $l$
such that $P_1,\,P_2,\,P_3$ are not at the same line and none of them
coincides with the intersection point $Q \in l_1 \cap l_2 \cap l_3$;
$l \in {\cal L}_b$ iff $P_1,\,P_2,\,P_3$ are at the same line, and
$l \in {\cal L}_c$ iff $P_i = Q$ for some $i$.
To prove the claim we may assume that $l \in {\cal L}_i,\,i=a,\,b,\,c,$ is one
of the
the standard
configurations $l_a,\,l_b,\,l_c$ described below. For all of them
$l_1=\{x=0\}$, $l_2=\{y=0\}$, $l_3=\{x=y\}$
in the homogeneous coordinates $(x:y:z)$ in $\pr^2$, and, respectively,

\ss

\noin $l_a$) $P_1=(0:1:0)$, $P_2=(1:0:0)$, $P_3=(1:1:1)$,

\noin $l_b$) $P_1=(0:1:0)$, $P_2=(1:0:0)$, $P_3=(1:1:0)$,

\noin $l_c$) $P_1=(0:0:1)$, $P_2=(1:0:0)$, $P_3=(1:1:0)$.

\ss

\noin  Let $C \in {\cal C}(d)$ be given by the equation
$$ \sum_{0\le i+j \leq d} a_{ij}x^iy^jz^{d-i-j} =0\,.$$ Then
$C \in {\cal C}(d)_l$ iff, respectively,

\ss

\noin a) $a_{0,d}=a_{0,d-1}=a_{0,d-2}=a_{d,0}=a_{d-1,0}=a_{d-2,0}=0$,

\noin  $\sum_{i+j \leq d} a_{i,j}=0$, $\sum_{i+j \leq d} (i+j)a_{i,j}=0$,
$\sum_{i+j \leq d} (i+j)(i+j-1)a_{i,j}=0$,

\noin b) $a_{0,d}=a_{0,d-1}=a_{0,d-2}=a_{d,0}=a_{d-1,0}=a_{d-2,0}=0$,

\noin $\sum_{i+j=d} a_{i,j}=0$, $\sum_{i+j=d-1} a_{i,j}=0$, $\sum_{i+j=d-2}
a_{i,j}=0$,

\noin c) $a_{0,0}=a_{0,1}=a_{0,2}=a_{d,0}=a_{d-1,0}=a_{d-2,0}=0$,

\noin $\sum_{i+j=d} a_{i,j}=0$, $\sum_{i+j=d-1} a_{i,j}=0$, $\sum_{i+j=d-2}
a_{i,j}=0$.

\ss

\noin Representing these equations on the Newton diagram, it is an easy
exercise to check that, if $d \geq 4$, they impose 9 independent conditions
on the coefficients $a_{ij}$ of $C$ in all the cases (a), (b) and (c),
and the claim follows. \qed

\vs

\noin To prove the proposition in the case when $d \ge 4$,
note that the subvariety
${\cal S}(d) := \bigcup_{l \in {\cal L}}  {\cal C}(d)_l \s {\cal C}(d)$
consists
of the orbits of the induced PGL$(3;\,\cz)$--action on ${\cal C}(d)$.
Moreover, it consists of the orbits of the subsets
${\cal C}(d)_{l_a},\,{\cal C}(d)_{l_b},\,{\cal C}(d)_{l_c}$. Since \\
dim$\,$PGL$(3;\,\cz) = 8$, due to the above claim,
all of these three orbits have codimension at least one. Hence,
the complement $ {\cal C}(d) \se {\cal S}(d)$ contains a Zariski open
subset. It remains to notice that the latter complement coincides with
the set of
smooth curves of degree $d$ whose inflection tangent lines are in general
position.

Consider further the remaining case $d = 3$. It is easily seen
that a plane cubic which satisfies one of the conditions (a) or (c)
is reducible. The only cubics which satisfy (b) are those from the linear
pencil $${\cal A} = \{C_{(\alpha : \beta)} = \{\alpha xy(x-y) - \beta z^3 =
0\},\,\,(\alpha : \beta)\in \pr^1\}\,.$$
Therefore,
${\cal C}(3)_{l_b} \s {\cal A}$. The linear pencil ${\cal A}$
is invariant under the action of the one parameter group of automorphisms
$(x : y : z) \longmapsto (x : y : cz),\,\,c \in \cz^*$. Hence, the
PGL$(3,\,\cz)$--orbit $S(3)$ of ${\cal C}(3)_{l_b}$ is of dimension
at most $8$. Once again, the complement ${\cal C}(3) \se S(3)$ is Zariski
open and it consists of the smooth cubics whose inflexional tangent lines
are in general position. This completes the proof. \qed

\vs

\begin{center} {\LARGE References} \end{center}

{\footnotesize
\noin [1] P. Aluffi, C. Faber. {\sl Linear orbits of $d$--tuples of points
in $\pr^1$}, J. reine angew. Math. 445 (1993), 205--220

\noin [2] A. B. Aure. {\sl Pl\"ucker conditions on plane rational curves},
Math. Scand. 55 (1984), 47--58, with {\sl Appendix} by S. A. Str\"omme,
ibid. 59--61

\noin [3] K. Azukawa, M. Suzuki. {\sl Some examples of algebraic degeneracy
and hyperbolic manifolds}, Rocky Mountain J. Math. 10 (1980), 655--659

\noin [4] J. A. Carlson, M. Green. {\sl Holomorphic curves in the plane}, 
Duke Math. J. 43 (1976), 1--9

\noin [5] A. I. Degtyarev. {\sl Topology of plane projective algebraic
curves}, PhD Thesis, Leningrad State University, 1987 (in Russian)

\noin [6] A. Degtyarev.  {\sl Quintics in ${\bf CP}^2$ with nonabelian
fundamental group}. Preprint, Max--Planck--Institut f\"ur Mathematik, Bonn,
1995, 21p.

\noin [7] P. Deligne. {\sl Le groupe fondamental du complement d'une
courbe plane n'ayant que des points doubles ordinaires est ab\'elien},
Sem. Bourbaki, 1979/1980, Lect. Notes in Math. vol. 842, Springer-Verlag
(1981), 1--25

\noin [8] G. Dethloff, G. Schumacher, P.-M. Wong. {\sl Hyperbolicity of the
complements of plane algebraic curves}, Amer. J. Math. 117 (1995), 573--599

\noin [9] G. Dethloff, G. Schumacher, P.-M. Wong.
{\sl On the hyperbolicity of the
complement of curves in algebraic surfaces: The three component case},
Duke Math. J. 78 (1995), 193--212

\noin [10] G. Dethloff, M. Zaidenberg. {\sl Examples of plane curves of low
degrees
with hyperbolic and C--hyperbolic complements.}
Proc. Conf. "Geometric Complex Analysis", Hayama, March 19--29, 1995
(to appear in World Scientific, Singapore)

\noin [11] I. Dolgachev, A. Libgober. {\sl On the fundamental group of the
complement to a discriminant variety}, In: Algebraic Geometry,
Lecture Notes in Math. 862, 1--25, N.Y. e.a.: Springer, 1981

\noin [12] H. Flenner, M. Zaidenberg. {\sl On a class of rational cuspidal
plane curves.} Preprint Mathematica G\"ottingensis 28 (1995), 1--31

\noin [13] W. Fulton. {\sl On the fundamental group of the complement to a
node curve.} Ann. of Math. (2) 111 (1980), 407-409

\noin [14] I. M. Gelfand, M. M. Kapranov, A.V. Zelevinsky.
{\sl Discriminants, resultants and multidimensional determinants},
Boston e.a.: Birkh\"auser, 1994

\noin [15] H. Grauert, U. Peternell. {\sl Hyperbolicity of the complements
of plane curves}, Manuscr. Math. 50 (1985), 429--441

\noin [16] M. Green. {\sl The complement of the dual of a plane curve and
some new hyperbolic manifolds}, in: 'Value Distribution Theory`,
Kujala and Vitter, eds., N.Y.: Marcel Dekker, 1974, 119--131

\noin [17] M. Green. {\sl Some examples and counterexamples in value
distribution theory}. Compos. Math. 30 (1975), 317-322

\noin [18] D. Hilbert. {\sl Theory of algebraic invariants}, Cambridge:
Cambridge Univ. Press, 1993

\noin [19] Sh. I. Kaliman. {\sl The holomorphic universal covers of spaces of
polynomials without multiple roots}, Selecta Mathem. form. Sovietica, 12 (1993)
No. 4, 395--405

\noin [20] J. Kaneko. {\sl On the fundamental group of the complement to
a maximal cuspidal plane curve}, Mem. Fac. Sci. Kyushu Univ. Ser. A. 39
(1985), 133-146

\noin [21] Sh. Kobayashi. {\sl Hyperbolic manifolds and holomorphic mappings},
N.Y. a.e.: Marcel Dekker, 1970

\noin [22] V. Ja. Lin. {\sl Liouville coverings of complex spaces, and
amenable groups}, Math. USSR Sbornik, 60 (1988), 197--216

\noin [23] K. Masuda, J. Noguchi. {\sl A construction of hyperbolic
hypersurface of $\pr^n(\cz)$}, Preprint 1994, 24 p. (to appear in
Mathem. Annalen)

\noin [24]  A. Nadel. {\sl Hyperbolic surfaces in $\pr^3$}, Duke Math.
J. 58 (1989), 749--771

\noin [25] M. Namba. {\sl Geometry of projective algebraic curves},
N.Y. a.e.: Marcel Dekker, 1984

\noin [26] M. Oka. {\sl Symmetric plane curves with nodes and cusps}, J. Math.
Soc. Japan, 44, No. 3 (1992), 375--414

\noin [27] A. Yu. Ol'shanskiy, A. L. Shmel'kin. {\it Infinite groups.}
In: Algebra IV, Encyclopaedia of Math. Sci. 37, Berlin e.a.: Springer, 1993,
3--95

\noin [28] R. Piene. {\sl Cuspidal projections of space curves}, Math. Ann.
256 (1981), 95--119

\noin [29] V. Popov, E. Vinberg. {\sl Invariant theory}, In: `Algebraic
Geometry IV', Parshin and Shafarevich, eds., 123-278, Berlin e.a.:
Springer, 1994

\noin [30] H. Royden. {\it Automorphisms and isometries of Teichm\"uller
space}, In: Advances in the Theory of Riemann Surfaces, 1969 Stony Brook
Conf. Ann. of Math. St. 66, Princeton, N.J.: Princeton Univ. Press, 1971,
369--383

\noin [31] F. Severi. {\sl Vorlesungen \"{u}ber algebraische Geometrie},
Leipzig: Teubner, 1921

\noin [32] Y.--T. Siu, S.--K. Yeung. {\sl Hyperbolicity of the complement
of a generic smooth curve of high degree in the complex projective plane},
preprint, 1994, 56 p.

\noin [33] G. Veronese. {\sl Behandlung der projectivischen Verh\"altnisse
der R\"aume von verschiedenen Dimensionen durch das Princip des Projicirens
und Schneidens}, Math. Ann. 19 (1882), 193--234

\noin [34] R. J. Walker. {\sl Algebraic curves}, Princeton Math.
Series 13, Princeton, N.J.: Princeton University Press, 1950

\noin [35] M. Zaidenberg. {\sl On hyperbolic embedding of complements of
divisors and the limiting behavior of the Kobayashi-Royden metric}, Math.
USSR Sbornik 55 (1986), 55--70

\noin [36] M. Zaidenberg. {\sl Stability of hyperbolic imbeddedness and
construction of examples}, Math. USSR Sbornik 63 (1989), 351--361

\noin [37] M. Zaidenberg. {\sl Hyperbolicity in projective spaces}, Proc.
Conf.  on Hyperbolic and Diophantine Analysis, RIMS, Kyoto, Oct. 26--30,
1992. Tokyo, TIT, 1992, 136--156

\noin [38] O. Zariski. {\sl Collected Papers}. Vol III : {\sl Topology of
curves and surfaces, and special topics in the theory of algebraic
varieties},  Cambridge, Massachusets e. a.: The MIT Press, 1978 \vs

\noindent Gerd Dethloff,\\
Mathematisches Institut der Universit\"at G\"ot\-tin\-gen,\\
Bunsenstrasse 3-5,\\
37073 G\"ot\-tin\-gen,\\
Germany.\\
e-mail: DETHLOFF@CFGAUSS.UNI-MATH.GWDG.DE

\vs

\noin  Mikhail Zaidenberg,\\
Universit\'{e} Grenoble I,\\
Institut Fourier des Math\'ematiques,\\
BP 74,\\
38402 St. Martin d'H\`{e}res--c\'edex,\\
France.\\
e-mail: ZAIDENBE@PUCCINI.UJF--GRENOBLE.FR}

\end{document}